\documentclass[a4paper,12pt]{article}
\usepackage[utf8]{inputenc}
\usepackage{amsmath}
\usepackage{amssymb}
\usepackage{amsthm}
\usepackage{xcolor}

\usepackage{enumerate}
\usepackage{amssymb, amsmath,graphicx}
\usepackage{todonotes}
\usepackage{verbatim}
\usepackage{float}

\theoremstyle{plain}
\newtheorem{thm}{Theorem}[section]
\newtheorem{theorem}[thm]{Theorem}
\newtheorem{proposition}[thm]{Proposition}
\newtheorem{definition}[thm]{Definition}
\newtheorem{lemma}[thm]{Lemma}
\newtheorem{corollary}[thm]{Corollary}

\newtheorem{conjecture}[thm]{Conjecture}
\newtheorem{open}[thm]{Question}

\newtheorem{remark}[thm]{Remark}

\newtheorem*{example*}{Example}
\newtheorem*{remark*}{Remark}

\newcommand{\Ent}{\mathrm{Ent}}

\newcommand{\Ch}{\mathrm{Ch}}

\newcommand{\Dom}{\mathcal{D}}

\newcommand{\Lip}{\mathrm{Lip}}

\newcommand{\cP}{\mathcal P}

\renewcommand{\d}{{\sf d}}

\newcommand{\W}{{\sf W}}
\newcommand{\mm}{\mathfrak m}

\def\N{{\mathbb N}}
\def\Z{{\mathbb Z}}
\def\R{{\mathbb R}}
\def\C{{\mathbb C}}

\def\Q{{\mathbb Q}}
\def\PP{{\mathbb P}}

\usepackage{amsmath, verbatim, amsfonts, amssymb, color}
\usepackage{mathrsfs}
\usepackage{dsfont}
 \usepackage{mathbbol}

\def\cF{\mathcal F}
\def\cE{\mathcal E}

\def\cP{\mathcal P}

\def\M{{\sf M}}%Riemannian manifold
\def\C{\mathcal{C}}%stetige Funktionen

\def\Lip{\mathrm{Lip}}%Lip.stetige Funktionen

\def\N{{\mathbb N}}    %Zahlkoerper
\def\Z{{\mathbb Z}} %Zahlkoerper
\def\R{{\mathbb R}}%Zahlkoerper
\def\Q{{\mathbb Q}}

\def\PP{\mathbb P}
\def\c{\mathfrak c}
\def\cc{^\c}
\def\Conj{\mathfrak C}
\def\K{\mathcal K}
\def\G{\mathcal G}
\def\Ch{{{\sf C\!h}}}
\begin{document}

\title{\bfseries Wasserstein Diffusion \\ on \\ Multidimensional Spaces% --- And Yet It Moves
}

\author{Karl-Theodor Sturm
\\[1cm]
\small Hausdorff Center for Mathematics \& Institute for Applied Mathematics\\
\small University of Bonn, Germany 
%sturm\a uni-bonn.de 
}

\maketitle

\abstract{Given any closed Riemannian manifold $\M$, we construct a reversible diffusion process on the space $\cP(\M)$ of probability measures on $\M$ that is
\begin{itemize}

\item reversible w.r.t.~the entropic measure $\PP^\beta$ on $\cP(\M)$, 
heuristically given as 
\begin{equation*}
d\mathbb{P}^\beta(\mu) =\frac{1}{Z} e^{-\beta \, \Ent(\mu \mid \mm)}\ d\mathbb{P}^*(\mu);
\end{equation*}
\item associated with a regular Dirichlet form 
with carr\'e du champ derived from the Wasserstein gradient in the sense of Otto calculus
$$\cE_\W(f)=\liminf_{g\to f}\ \frac12\int_{\cP(\M)} \big\|\nabla_\W g\big\|^2(\mu)\ d\PP^\beta(\mu);$$
\item non-degenerate, at least in the case of the $n$-sphere and the $n$-torus.
\end{itemize}}

\bigskip\bigskip

{\hfill {\it And yet it moves.}}

{\hfill \tiny{Galileo Galilei, 1633} }

 \section{Introduction}
 \paragraph{(a)}
 Given any closed Riemannian manifold $\M$ with normalized volume measure $\mm$, we construct a reversible diffusion process on the space $\cP(\M)$ of probability measures on $\M$ that is
 reversible w.r.t.~the entropic measure $\PP^\beta$ on $\cP(\M)$, a random measure 
heuristically given as 
\begin{equation}\label{formal}
d\mathbb{P}^\beta(\mu) =\frac{1}{Z} e^{-\beta \, \Ent(\mu \mid \mm)}\ d\mathbb{P}^*(\mu)
\end{equation}
in terms of the relative entropy $\Ent(. \mid \mm)$ and 
a (non-existing) uniform distribution $\PP^*$ on $\cP(\M)$. This ansatz will be made rigorous by means of the conjugation map, a continuous involutive bijection 
$$\Conj_\cP: \ \begin{array}{ r c c c }
& \cP(\M)&\to& \cP(\M)\\
&\exp(\nabla\varphi)_*\mm&\mapsto &\exp(\nabla\varphi^\c)_*\mm,\\
\end{array}$$
the definition of which is based on the fact that, by the McCann-Brenier Theorem, every measure $\mu\in \cP(\M)$ can be represented as $\mu=\exp(\nabla\varphi)_*\mm$ in terms of a $\c$-convex function $\varphi$. 
The conjugate measure %$\mu^\c:=\Conj(\mu)$ 
is then represented %as $\mu^\c=\exp(\nabla\varphi^\c)_*\mm$ 
in terms of
 $\varphi^\c$, the $\c$-conjugate function for the given $\varphi$ and the cost function $\c(x,y):=\frac12\d^2(x,y)$.
 
 Rigorously, the entropic measure is then defined as
the push forward of the Dirichlet-Ferguson measure $\mathbb{Q}^\beta$ (with reference measure $\beta \mm$) under the conjugation map
 $$ \mathbb{P}^\beta := (\Conj_\mathcal{P})_\ast \mathbb{Q}^\beta.$$
 This way, our heuristic ansatz \eqref{formal} will finally result in a Large Deviation Principle with rate function 
 $\Ent(\, .\mid\mm)$. % for  the entropic measures $\PP^\beta$.
 
  \begin{theorem} For every $A\subset\cP(\M)$,
 \begin{equation*}
-\inf_{\mu\in A^\circ}\Ent(\mu\mid\mm)\le\liminf_{\beta\to\infty}\frac1\beta\log\PP^\beta(A)\le \limsup_{\beta\to\infty}\frac1\beta\log\PP^\beta(A)\le -\inf_{\mu\in \bar A}\Ent(\mu\mid\mm).
\end{equation*}
\end{theorem}
Despite the regularizing effect of the relative entropy, the entropic measure 
gives mass only to singular measures.

 \begin{theorem} Assume that $\M$ is the $n$-sphere, the $n$-torus or a convex Euclidean compactum. 
% Then for every $\beta>0$ the entropic measure $ \mathbb{P}^\beta$
%\begin{itemize}
%%\item has full topological support on $ \cP(\M)$;
%\item gives mass only to measures $\mu$ on $\M$ which have no absolutely continuous part;
%\item  gives mass only to measures %$\mu$ on $\M$ 
%which have no atoms.
%\end{itemize}
%
 Then  $ \mathbb{P}^\beta$-a.e.~$\mu$ has
%\begin{itemize}
%\item has full topological support on $ \cP(\M)$;
  no absolutely continuous part and
  no atoms.
%\end{itemize}
\end{theorem}

\medskip

\begin{center}
\includegraphics[width=4.cm]{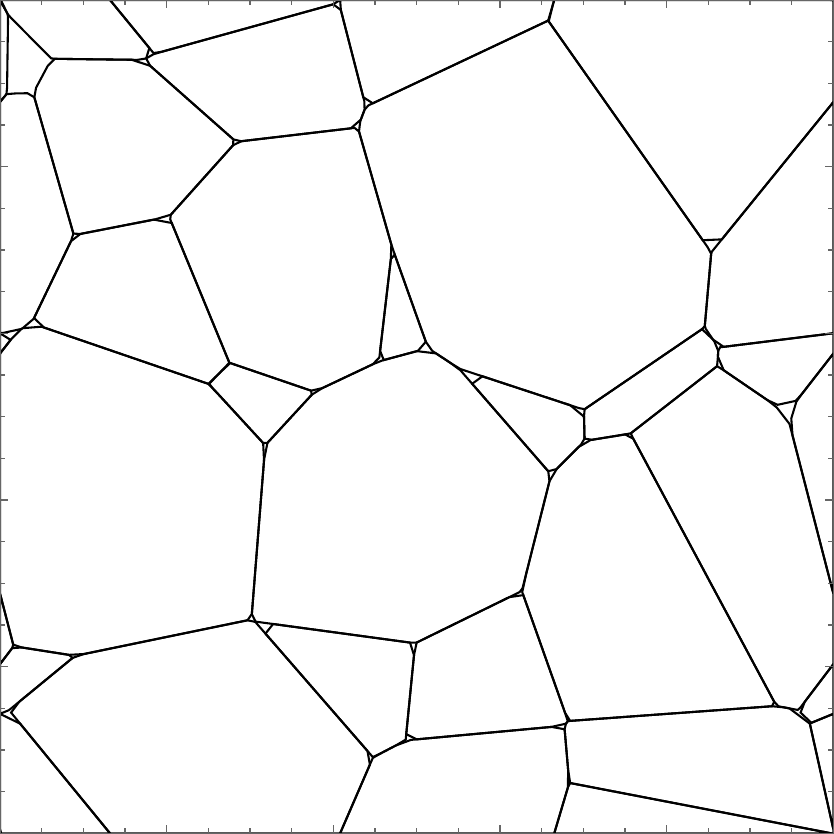}\qquad
\includegraphics[width=4.cm]{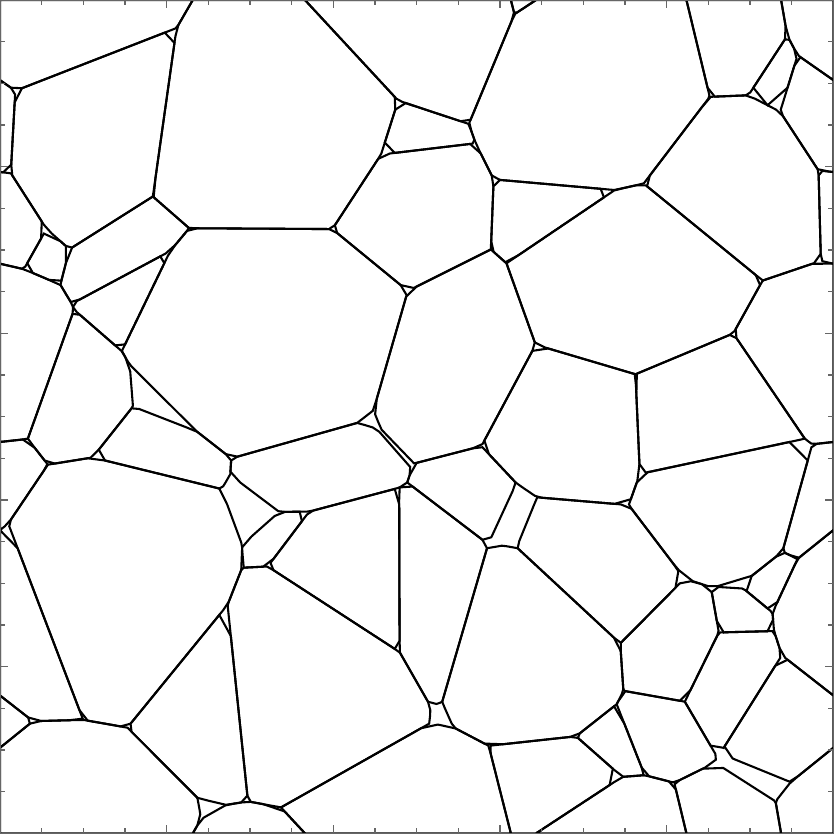}\qquad
\includegraphics[width=4.cm]{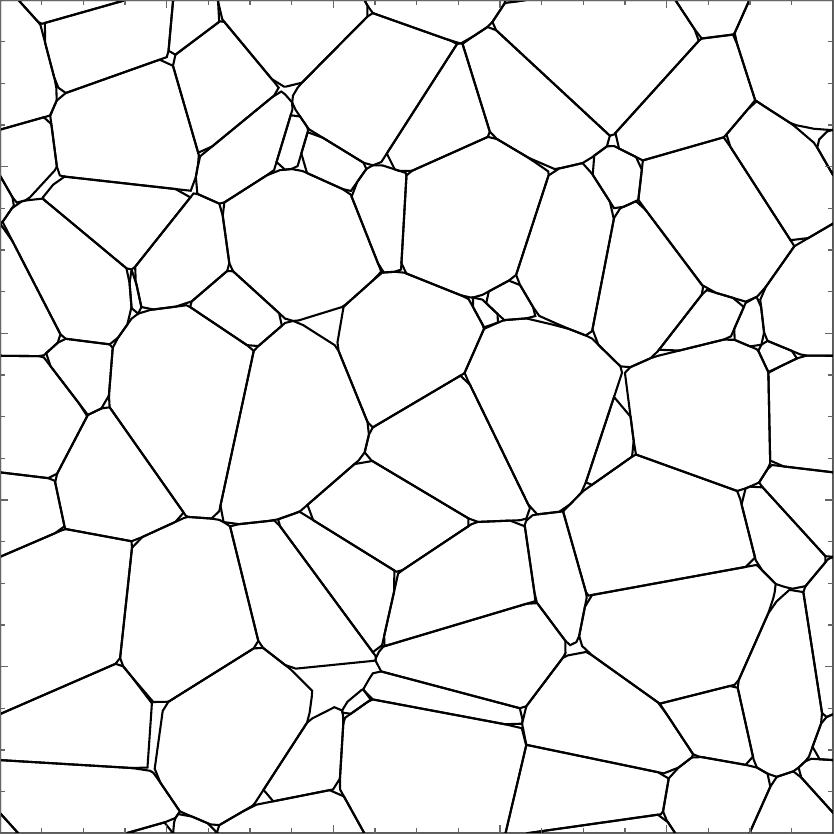}

\medskip

{{Figure 1:} Support %of realizations 
of $\PP^\beta$-distributed $\mu\in\cP(\R^2/\Z^2)$ for $\beta=15$, $30$ and $75$, resp.,\\
simulations provided by courtesy of Lorenzo Dello Schiavo.}
\end{center}
%\end{figure}

%A typical realization of $\mu$ looks like the simulations of the early universe in \cite{early}.

%\begin{figure}
\paragraph{(b)}
Our next goal is to introduce a `canonical' \emph{Wasserstein Dirichlet form} $\cE_\W$ and associated \emph{Wasserstein diffusion process} on $\cP(\M)$. We will define the former as the relaxation on $L^2(\cP(\M), \PP^\beta)$  of the pre-Dirichlet form 
\begin{align*}
{\mathcal E}^0_\W(f):=\frac12\int_{{\mathcal P}(\M)}
 \big\|\nabla_\W f\big\|^2(\mu)
 \, d\PP^\beta(\mu),
\end{align*}
defined on {cylinder functions on  $\cP(\M)$},
where $\nabla_\W$ denotes the Wasserstein gradient  in the sense of the Otto calculus and where $\PP^\beta$ is the entropic measure with parameter $\beta>0$.

Recall that a function $f:  {\mathcal P}(\M)\to\R$ is called \emph{cylinder function} if
$$
f(\mu)=F\left(\int_\M \vec V\,d\mu\right)$$
for suitable $ k\in\N, F\in\C^1(\R^k)$ and $\vec V=(V_1,\ldots,V_k)\in \C^1(\M,\R^k)$. 
The squared norm of the \emph{Otto-Wasserstein gradient} at the point $\mu\in\cP(\M)$ is given
 is given by 
$$\big\|\nabla_\W f\big\|^2(\mu):=
\sum_{i,j=1}^k \Big(\partial_i F\cdot \partial_j F\Big)\bigg(\int_\M \vec V\,d\mu\bigg)
\int_\M \langle\nabla V_i,\nabla V_j\rangle\,d\mu
.$$

\begin{theorem}\begin{enumerate}[\rm (i)]
\item  $\cE_\W$ is a strongly local, regular Dirichlet form on $L^2(\cP(\M), \PP^\beta)$.

\item It has intrinsic metric 
${\sf d}_\W\ge \W_2$ and admits a  carr\'e du champ $\Gamma_\W$. For Lipschitz functions the latter is dominated by the squared asymptotic Lipschitz constant.

\item There exists  a reversible strong Markov process $\big((\rho_t)_{t\ge0}, ({\mathbf P}_\mu)_{\mu\in\cP(\M)}\big)$ with a.s.~continuous trajectories which is 
properly associated with ${\mathcal E}_\W$.
\end{enumerate}
\end{theorem}
This stochastic process is called \emph{Wasserstein diffusion.}

\paragraph{(c)}
One of the main results now is that the Wasserstein diffusion is non-degenerate, that is, it is really moving:
$$\PP^\beta\Big\{\mu: \ {\mathbf P}_\mu\big\{\exists t>0: \rho_t\not= \rho_0\big\}>0\Big\}>0,$$
or equivalently, $\cE_\W\not\equiv0$ on its domain --- at least in the case of the sphere  and the torus.

\begin{theorem} On the $n$-sphere and the $n$-torus, the Wasserstein Dirichlet form is non-vanishing: 
$$\mathcal E_\W(f)\not=0$$ for every non-vanishing
$f:{\mathcal P}(\M)\to\R$ which is
\begin{enumerate}[\rm (i)]
\item antisymmetric (i.e. $f\circ \hat\Phi_1=-f$ for suitable Lipschitz isometries, % on $\M$, 
see Def.~\ref{def-anti}) or

\item given as
$f(\mu):= \int_\M V\,d\mu$
in terms of a Lipschitz function $V:\M\to\R$.
\end{enumerate}
\end{theorem}

\section{Conjugation Map}
\subsection{$\c$-Convexity}
Throughout the sequel, $\M$ will be a compact subset of a complete Riemannian manifold $\hat\M$ with Riemannian distance $\d$. 
Furthermore, $\mm$ will denote a probability measure with full topological support $\M$ which is absolutely continuous with respect to the volume measure and satisfies a Poincar\'e inequality: $\exists c > 0$
$$\int |\nabla u|^2 d\mm \ge  c \int u^2\d\mm$$
for all weakly differentiable $u : \M\to \R$ with $\int ud\mm = 0$.
For closed Riemannian manifolds, there is a canonical choice for m, namely, the normalized Riemannian volume measure.

Such a set $\M$ will be called \emph{Riemannian compactum} and such a measure $\mm$ \emph{reference measure}.
The \emph{cost function} will always be
$$\c(x,y):=\frac12 \d^2(x,y).$$

\begin{definition}
A function $\varphi: \M \rightarrow \R $ is called $\c$-convex if there exists a function $\psi: \M \rightarrow \R $ such that
\begin{displaymath}
 \varphi(x)=-\underset{y\in M}{\inf} \left[\frac{1}{2} \d^2(x,y)+\psi(y)\right]
\end{displaymath}
for all $x \in \M$. In this case, $\varphi$ is called %\emph{generalized Legendre transform} of $\psi$ or
 \emph{conjugate} of $\psi$ and denoted by
$$\varphi=\psi\cc.$$
\end{definition}
In our presentation, the ${}^\c$ stands both for the quadratic cost function $\c$ as above and for `conjugate'. 

 A function $\varphi$ is $\c$-convex if and only if the function $-\varphi$ is $\c$-concave in the sense of
\cite{Rockafellar70, Rueschendorf96, McCann01,Villani08}. A function $\varphi$ on $\R^n$ is $\c$-convex if and only if the function $x\mapsto \varphi(x)+\frac12|x|^2$ is convex in the usual sense.
Let us summarize some of the basic facts on $\c$-convex functions. See \cite{Rockafellar70}, \cite{Rueschendorf96}, \cite{McCann01} and \cite{Villani08} for
details.

\begin{lemma}
\begin{enumerate}[\rm (i)]
 \item A function $\varphi$ is $\c$-convex if and only if ${\varphi\cc}\cc = \varphi$.
 \item Every $\c$-convex function is bounded, Lipschitz continuous and differentiable almost everywhere with gradient bounded by $D= \underset{x,y \in
     M}{\sup}
     d(x,y)$.
\end{enumerate}
\end{lemma}

In the sequel, ${\K}={\K}(\M)$ will denote the set of $\c$-convex functions on $\M$ and $\tilde{\K} = \tilde{\K}(\M)$ will denote the set of equivalence classes
in  ${\K}$ with $\varphi_1 \sim \varphi_2$ iff $\varphi_1 - \varphi_2$ is constant.
${\K}$ will be regarded as a subset of the Sobolev space $H^1=H^1(\M,\mm)$ with norm
$\| u \|_{H^1} = \left[ \int_\M  |\nabla u|^2 d\mm + \int_\M u^2 d\mm \right]^{1/2}$
and $\tilde{\K} = {\K}/const $ will be regarded as a subset of the space $\tilde{H^1}= H^1 / const$ with norm
$\| u \|_{\tilde{H^1}} = \left[ \int_\M | \nabla u|^2 d\mm \right]^{1/2}$.

\begin{lemma}
For each Borel map $g:\M\rightarrow \M$ the following are equivalent:
\begin{enumerate}[\rm (i)]
 \item $\exists \varphi \in \tilde\K : g= \exp(\nabla \varphi)$ \ $\mm$-a.e.~on $\M$;
 \item $g$ is an optimal transport map from $\mm$ to $g_\ast \mm$ in the sense that it is a minimizer of $h \mapsto \int_\M \d^2(x,h(x))\mm(dx)$ among all Borel maps
     $h:\M\rightarrow \M$ with $h_\ast \mm = g_\ast \mm$.
\end{enumerate}
\end{lemma}
In this case, the function $\varphi \in \tilde \K$ in (i) is defined uniquely. Moreover, in (ii) the map $g$ is the unique minimizer of the given minimization
problem.

A Borel map $g:\M\rightarrow \M$ satisfying the properties of the previous proposition will be called \emph{monotone map} or \emph{optimal Lebesque transport}.
The set of $\mm$-equivalence classes of such maps will be denoted by $\G=\G(\M)$ and regarded as a subset of the space of maps $L^2((\M,\mm);(\M,\d))$ with metric $\d_2(f,g) = \left[ \int_\M \d^2(f(x),g(x))\mm(dx) \right]^{1/2} $.
Note that $\G(\M)$ does not depend on the choice of $\mm$ (as long as m is absolutely continuous with full support)!

According to our definitions, the map $\Upsilon:\varphi \mapsto \exp(\nabla \varphi)$ defines a bijection between $\tilde \K$ and $\G$.
Recall that $\mathcal{P} = \mathcal{P}(\M)$ denotes the set of probability measures $\mu$ on $\M$ (equipped with its Borel $\sigma$-field).

\begin{proposition}
 The map $\chi: g \mapsto g_\ast \mm$ defines a bijection between $\G$ and $\mathcal{P}(\M)$. That is, for each $\mu \in \mathcal{P}$ there exists a unique $g \in \G$
 ---
 called \emph{Brenier map} of $\mu$ --- with $\mu = g_\ast \mm$.
\end{proposition}

Due to the previous observations, there exist canonical bijections $\Upsilon$ and $\chi$ between the sets $\tilde \K$, $\G$ and $\mathcal{P}$.
Actually, these bijections are even homeomorphisms with respect to the natural topologies on these spaces.
\begin{lemma}[{ \cite[Prop.~2.5]{Sturm-multi}}]\label{prop2.5}
Consider any sequence $\left\lbrace \varphi_n \right\rbrace_{n \in \mathbb{N}} $ in $\tilde \K$ with corresponding sequences $\left\lbrace g_n \right\rbrace_{n \in
\mathbb{N}} = \left\lbrace \Upsilon(\varphi_n) \right\rbrace_{n \in \mathbb{N}}$ in $\G$ and $\left\lbrace \mu_n \right\rbrace_{n \in \mathbb{N}} = \left\lbrace
\chi(g_n) \right\rbrace_{n \in \mathbb{N}}$ in $\mathcal{P}$ and let $\varphi \in \tilde \K$, $g=\Upsilon(\varphi) \in \G$, $\mu =\chi(g) \in \mathcal{P}$. Then the
following are equivalent:
\begin{enumerate}[\rm (i)]
 \item $\varphi_n \longrightarrow \varphi$ in $\tilde H_1$
 \item $g_n \longrightarrow  g$ in $L^2((M,\mm);(M,d))$
 \item $g_n \longrightarrow  g$ in $\mm$-probability on $M$
 \item $\mu_n \longrightarrow  \mu$ in $L^2$--Kantorovich--Wasserstein distance $\mathsf W_2$
 \item $\mu_n \longrightarrow  \mu$ weakly.
\end{enumerate}

\end{lemma}

\subsection{Conjugation Map}

Let us now consider the conjugation map $ \Conj_\K : \varphi \mapsto \varphi \cc$ in more details. It acts on functions $\varphi : \M\rightarrow \R$ and is defined as follows
\begin{equation*}
 \varphi\cc (x) = -\underset{y\in M}{\inf} \left[ \frac{1}{2} \d^2(x,y) + \varphi(y) \right] .
\end{equation*}
Obviously, $\Conj_\K$ maps $\K$ bijectively onto itself and
$\Conj^2 _\K = Id$. Moreover,  $\Conj_\K (\varphi + \lambda) = \Conj_\K (\varphi) - \lambda$ for each $\lambda \in \R$.
Hence, $\Conj _\K$ extends to a bijection $\Conj_{\tilde \K}: \tilde \K \rightarrow \tilde \K$.
Composing this map with the bijections $\chi:\G \rightarrow \mathcal{P}$ and $\Upsilon : \tilde \K \rightarrow \G$ we obtain involutive bijections
\begin{equation*}
 \Conj_\G = \Upsilon \circ \Conj_{\tilde \K} \circ \Upsilon^{-1}:\G \rightarrow \G, \qquad\quad
  \Conj_\mathcal{P} = \chi \circ \Conj_{\G} \circ \chi^{-1}:\mathcal{P} \rightarrow \mathcal{P},
\end{equation*}
called \emph{conjugation map} on $\G$ or on $\mathcal{P}$, respectively.
(Note that the maps $\chi$ and $ \Conj_\mathcal{P} $  strongly depend on the choice of the measure $\mm$.)
Given a monotone map $g \in \G$, the monotone map
\begin{equation*}
 g\cc := \Conj_\G (g)
\end{equation*}
will be called \emph{conjugate map} or \emph{generalized inverse map}; given a probability measure $\mu \in \mathcal{P}$ the probability measure
\begin{equation*}
 \mu\cc := \Conj_\mathcal{P} (\mu)
\end{equation*}
will be called \emph{conjugate measure}.

\begin{lemma}
 The conjugation map $$\Conj_\K : \K \rightarrow \K$$ is continuous.
\end{lemma}

This result was already stated  as Lemma 3.5 in \cite{Sturm-multi}. However, the proof there is incomplete as kindly pointed out by  Mathav Murugan.

\begin{proof}
To simplify notation denote $\Conj_\mathcal{K}$ by $\Conj$.

(i) Let us consider a sequence $(\varphi_\ell)_{\ell \in \mathbb{N}}$ in $\K$ with $\varphi_\ell \rightarrow \varphi$ in $H^1(M)$.  Outside of some $\mm$-zero set $M_0\subset M$, the convergence $\varphi_\ell \rightarrow \varphi$ will also be pointwise.
 Choose a countable  set $\left\lbrace y_i\right\rbrace_{i \in \mathbb{N}}$ in $M\setminus M_0$, dense in $M$, and for $j \in \mathbb{N}$ define $\Conj_j(\psi): M\to\mathbb{R}$ for $\psi \in \K$  by $\Conj_j(\psi)(x) = - \underset{i=1,\ldots,j}{\inf} [\frac{1}{2} d^2(x,y_i)+\psi(y_i) ]$.
Then for each $j \in \mathbb{N}$ as
$\ell\rightarrow\infty$
$$\Conj_j (\varphi_\ell) \rightarrow \Conj_j (\varphi)$$ pointwise on $M$.
%
%Moreover, as $j\rightarrow\infty$ $${\psi\cc}^j \nearrow \psi\cc \qquad \text{ pointwise on M} .$$
For each $\varepsilon > 0$ choose $j=j(\varepsilon) \in \mathbb{N}$ such that the set $\left\lbrace y_i\right\rbrace _{i=1,\ldots,j(\varepsilon)}$ is an
$\varepsilon$-covering of the compact space $M$. Then
\begin{align*}
\mid \Conj_j(\psi)(x) - \Conj(\psi)(x) \mid  \, &\leq
\sup_{y \in M}\, \inf_{i=1,\ldots,j} \Big| \frac{1}{2} d^2 (x,y) - \frac{1}{2} d^2(x,y_i) + \psi(y) - \psi(y_i) \Big| \, \\
&\leq \sup_{y \in M}\, \inf_{i=1,\ldots,j} 2 D\cdot d(y,y_i) \leq 2 D \varepsilon\end{align*}
uniformly in $x \in M$ and $\psi \in\K$.
Recall that each $\psi \in \K$ is Lipschitz continuous with Lipschitz constant $D$.
The uniform convergence of $\Conj_j \rightarrow \Conj$ together with the previous pointwise convergence 
$\Conj_j (\varphi_\ell) \rightarrow \Conj_j (\varphi)$ 
 implies $$\Conj(\varphi_\ell) \rightarrow \Conj(\varphi)$$  as
$\ell\rightarrow\infty$ pointwise on $M$  and thus also in $L^2(M)$. 
%In particular, 
%$\Phi_\ell=\Conj(\varphi_\ell)-\int_M\Conj(\varphi_\ell)d\mm$ converges to $\Phi=\Conj(\varphi)-\int_M\Conj(\varphi)d\mm$ in $L^2(M)$.

\medskip

(ii) It remains to prove that $\Conj(\varphi_\ell) \rightarrow \Conj(\varphi)$ in $H^1$ as $\ell\to\infty$.
For this purpose, we choose a constant $A$ and cover $M$ by finitely many open $U_1,\ldots, U_k$ such that for each $j=1,\ldots,k$,
\begin{itemize}
\item
$\text{Hess}_x d^2(x,y)\le A\,\text{Id} \ $  for all $x,y\in M$
\item
$ \text{Hess}_x d^2(x,z)\ge  \text{Id} \ $  for all $x,z\in U_j$
\item there exists geodesic vector fields $X_1,\ldots,X_n$ on $U_j$ which at each $x\in U_j$ span the tangent space $T_xM$.
\end{itemize}
Let us fix $j=1,\ldots,k$ and a (`reference') point $z\in U_j$.
For $\psi\in\K$, define 
$$\bar\psi(x):=\psi(x)+A\,d^2(x,z).$$
Then $\bar\psi$ is \emph{convex} on $U_j$.
Indeed, since $\psi(x)= \underset{y\in M}{\sup} \big[-\frac{1}{2} d^2(x,y)-\psi^\c(y) \big]$ and due to our assumptions,
\begin{align*}
 \text{Hess}_x\bar\psi(x)\ge \inf_{y\in M}  \text{Hess}_x \Big[-\frac{1}{2} d^2(x,y)\Big]+A\,\text{Hess}_xd^2(x,z)\ge0
\end{align*}
for every $x\in U_j$.

\medskip

(iii) To simplify notation, we write  $\Phi_\ell:=\Conj_j (\varphi_\ell)+ A\,d^2(\,.\,,z)$ and $\Phi:=\Conj_j (\varphi)+ A\,d^2(\,.\,,z)$ for
given $\varphi_\ell$  and $\varphi$ as above.
 It suffices to prove that $\Phi_\ell\to\Phi$ in $H^1(U_j)$ as $\ell\to\infty$ where all the $\Phi_\ell$ and  $\Phi$ are convex and $D$-Lipschitz continuous functions on $U_j$ with
$\Phi_\ell\to\Phi$ pointwise as $\ell\to\infty$.
The convergence $\Phi_\ell\to\Phi$ in $H^1(U_j)$ will follow from the convergence $X_i\Phi_\ell\to X_i\Phi$ in $L^2$ for each $i=1,\ldots,n$.

To deduce the latter, let us fix $i$.
The geodesic vector field $X_i$ provides a decomposition of $U_j$ into geodesics. More precisely, there exists an open set $\Omega\subset \R^{n-1}$, open intervals $I_\omega$, and unit speed geodesics $\gamma_\omega: I_\omega\to U_j$ 
%($\forall \omega\in\Omega$) 
such that
$$U_j= \dot{\bigcup_{\omega\in\Omega}}\gamma_\omega(I_\omega).$$
Without restriction, $I_\omega= (-\delta,\delta)=:I$ for all $\omega$ and some $\delta>0$. Thus, $U_j$ can be identified with the set  $\Omega\times (-\delta,\delta)\subset \R^n$. %and the Riemannian volume measure on $U_j$

\medskip

(iv) 
For any $\omega\in \Omega$ now consider the functions
$$\Phi_\omega:=\Phi\circ \gamma_\omega: I\to U_j, \qquad \Phi_{\ell,\omega}:=\Phi_\ell\circ \gamma_\omega: I\to U_j.$$
Each of them is convex and $D$-Lipschitz continuous.
For $x\in U_j$, parametrized as $x=(\omega,t)\in \Omega\times (-\delta,\delta)$,
$$X_j\Phi(x)= \frac{\partial}{\partial t}\Phi_\omega(t)=:\Phi_\omega'(t), \qquad
X_i\Phi_\ell(x)= \frac{\partial}{\partial t}\Phi_{\ell,\omega}(t)=:\Phi_{\ell,\omega}'(t).$$
Pointwise convergence $\Phi_{\ell,\omega}(t)\to \Phi_{\omega}(t)$ as $\ell\to\infty$ and convexity of both $\Phi_{\ell,\omega}$ and $\Phi_{\omega}$ thus imply that for every $\omega$ and every $t$ under consideration, 
\begin{align*}
\limsup_{\ell\to\infty}\Phi'_{\ell,\omega}(t)&\le\limsup_{\ell\to\infty} \frac1\epsilon\Big[\Phi_{\ell,\omega}(t+\epsilon)-\Phi_{\ell,\omega}(t)
\Big]\\
&= \frac1\epsilon\Big[\Phi_{\omega}(t+\epsilon)-\Phi_{\omega}(t)\Big]\le \Phi'_{\omega}(t+\epsilon)
\end{align*}
for every $\epsilon>0$, and analogously,
\begin{align*}
\liminf_{\ell\to\infty}\Phi'_{\ell,\omega}(t)&\ge \Phi'_{\omega}(t-\epsilon).
\end{align*}
Since $\lim_{\epsilon\to0}\Phi'_{\omega}(t+\epsilon)-\Phi'_{\omega}(t-\epsilon)=0$ for $\mathfrak L^1$-a.e. $t\in I$ by convexity of $\Phi_\omega$, it follows that for $\mathfrak L^1$-a.e. $t\in I$,
$$\lim_{\ell\to\infty}\Phi'_{\ell,\omega}(t)\text{ exists and coincides with }\Phi'_{\omega}(t).$$
In other words, for every $\omega\in\Omega$ and $\mathfrak L^1$-a.e. $t\in I$,
$$X_i\Phi_\ell(\omega,t)\to X_i\Phi(\omega,t) \ \text{as }\ell\to\infty.$$
Thus, in particular, for $\mm$-a.e. $x\in U_j$,
$$X_i\Phi_\ell(x)\to X_i\Phi(x) \ \text{as }\ell\to\infty$$
which in turn implies the requested $L^2$-convergence since all functions under consideration are bounded by $D$.
\end{proof}

\begin{theorem}
 The conjugation map $$\Conj_\mathcal{P} :\mathcal{P}\rightarrow \mathcal{P}$$ is continuous (with respect to the weak topology).
\end{theorem}

\begin{proof} Let us first prove continuity of the conjugation map $\Conj_{\tilde \K}: \tilde \K \rightarrow \tilde \K$ (with respect to the $\tilde
H^1$-norm on $\tilde \K$).
Indeed, this follows from the previous continuity result together with the facts that the embedding $H^1\to\tilde H^1,\ \varphi\mapsto\tilde\varphi=\{\varphi+c:
c\in\mathbb{R}\}$ is continuous (trivial fact) and that the map $\tilde H^1\to H^1,\ \tilde\varphi=\{\varphi+c: c\in\mathbb{R}\}\mapsto\varphi-\int_M\varphi d\mm$ is
continuous (consequence of Poincar\'e inequality).

This in turn implies, due to Lemma \ref{prop2.5}, 
%that the conjugation map $\Conj _\G : \G \rightarrow \G $ is continuous (with respect to the $L^2$-metric on $\G$).
%Moreover, due to the same Proposition it therefore also implies 
that the conjugation map $$\Conj_\mathcal{P} :\mathcal{P}\rightarrow \mathcal{P}$$ is continuous (with
respect to the weak topology).
\end{proof}

%\begin{theorem}[{\cite[Thm.~3.6]{Sturm-multi}}]
% The conjugation map $$\Conj_\mathcal{P} :\mathcal{P}\rightarrow \mathcal{P}$$ is continuous (with respect to the weak topology).
%\end{theorem}

\begin{remark}\rm
As a corollary, we obtain that  the  conjugation map %s $\Conj_{\tilde \K}: \tilde \K \rightarrow \tilde \K$ and 
$\Conj_\mathcal{G} :\mathcal{G}\rightarrow \mathcal{G}$ is continuous as well.

In dimension $n=1$, the conjugation map $\Conj_\G$ is even an isometry from $\G$, equipped with the $L^1$-metric, into itself.
\end{remark}

\subsection{Conjugate Measures}

\begin{lemma} Assume that $\M$ is the $n$-sphere or the $n$-torus. 
Then for every $z\in\M$ there exists a unique maximizer $z^\c$ of $x\mapsto \d(x,z)$. The pair $z,z^\c$ is conjugate in the sense that it satisfies for all $y\in\M$,
\begin{equation}
\sup_{x\in\M}\Big[\d^2(x,z)-\d^2(x,y)\Big]= \d^2(z,z^\c)-\d^2(y,z^\c).
\end{equation}  
\end{lemma} 

\begin{proof}
In the case of the sphere, each of the points $x$ and $y$ lies on a geodesic from $z$ to $z^\c$. For given $y$, the quantity
\begin{equation}\label{3d2}
\d^2(x,z)+\d^2(y,z^\c)-\d^2(x,y)
\end{equation}
will be maximized among all $x$ with fixed distance from $z$ if $x$ lies on the same geodesic as $y$. If $x$  lies on the same geodesic as  $y$, then the above quantity \eqref{3d2} will be maximized for $x=z^\c$.

\medskip

In the case of the torus, we may assume that $z=(0,\ldots,0)\in (-\frac12,\frac12]^n$, that $z^\c=(\frac12,\ldots,\frac12)$, and that $y\in (0,\frac12]^n$. To maximize the quantity \eqref{3d2}, the point $x$ has to lie in the same subcube as $y$, i.e. $x\in (0,\frac12]^n$. But then all distance under consideration are just the distances of the ambient $\R^n$. Thus in this case,
\begin{align*}
|x-z|^2+|y-z^\c|^2-|x-y|^2-|z-z^\c|^2=2(x-z^\c)\cdot (y-z)\le0
\end{align*}
and of course equality holds for the choice $x=z^\c$.
\end{proof}

\begin{proposition}\label{conj-Diracs} Assume that $\M$ is the $n$-sphere or the $n$-torus. Then for every $z\in\M$,
\begin{equation}
(\delta_z)^\c=\delta_{z^\c}.
\end{equation}
\end{proposition}

\begin{proof} For every $z\in\M$, the function $\varphi=-\d^2(.,z)/2$ is the Brenier potential for the Dirac mass at $z$, i.e.
$$\delta_z=\exp\Big(-\nabla \d^2(.,z)/2\big)_*m.$$
According to the previous lemma,
$$\varphi^\c=-\d^2(.,z^\c)/2+ const.$$
Thus
$$(\delta_z)^\c=\exp\Big(-\nabla \d^2(.,z^\c)/2\big)_*m=\delta_{z^\c}.$$
\end{proof}

\begin{lemma}\label{densities}
 Let $\mu =g_\ast \mm \in \mathcal{P}$ and $\nu  = \mu\cc=
 f_\ast \mm$ with $g= \exp(\nabla\varphi)$ and $f= \exp(\nabla\varphi\cc)$ for some $\c$-convex $\varphi$. Assume that $\mu$ is absolutely continuous with density $\eta = \frac{d\mu}{d\mm}>0$ a.s.
  \begin{enumerate}[\rm (i)]

\item Then $$f(g(x)) = g(f(x)) =x \qquad \text{ for $\mm$-a.e.~} x\in \M.$$
\item Furthermore, $\nu$ is absolutely continuous with  density $\rho = \frac{d\nu}{d\mm}>0$ a.s.~and
\begin{equation*}
 \eta(x)= \frac1{\rho(f(x))}, \qquad   \rho(x)=\frac1{\eta(g(x))} \qquad \text{ for $\mm$-a.e.~} x \in \M.
\end{equation*}
%\item
%If $\mu$ and $\nu$ are absolutely continuous, then the Jacobian $\det Df(x)$ and $\det Dg(x)$ exist for almost every $x \in \M$ and satisfy
%\begin{equation*}
% \det D f(g(x)) \cdot \det D g(x) = \det D f(x) \cdot \det D g(f(x)) = 1,
%\end{equation*}
%\begin{equation*}
%\eta(x) = \det Df(x),
%\qquad
%\rho(x) = \det Dg(x)
%\end{equation*}
%for almost every $x \in \M$.
\end{enumerate}\end{lemma}

\begin{proof}  (i) 
 Absolute continuity of $\mu$ implies that there exists an optimal transport map which forwards it to $\mm$. Indeed, this will be the map $f= \exp(\nabla\varphi)$. In particular, $\varphi$ will be differentiable $\mu$-a.e. Positivity of the Radon-Nikodym density implies that the differentiability of $\varphi$ holds $\mm$-a.e.
 
 By optimality of $g$ and $f$ for the transports $\mm\to\mu$ and $\mu\to\mm$, they will also be optimal for any sub-transport of these transports. Thus for every measurable $A\subset \M$, 
 $$f(g(x))\in A \text{ for }\mm\text{-a.e. }x\in A.$$ 
 Hence, $f(g(x))=x$ for $\mm$-a.e.~$x\in \M$.
 Interchanging the roles of $f$ and $g$ yields the second identity.
 
 (ii)  [{\cite[Prop.~3.3]{Sturm-multi}}].
\end{proof}

\subsection{Entropy of Conjugate Measures}

\begin{definition}\label{def-ent} The \emph{relative entropy} of a probability measure $\mu\in\cP(\M)$ w.r.t.~the reference measure $\mm$ is given by
$$\Ent(\mu\mid\mm):=\left\{
\begin{array}{ll}
\int_\M \eta\log \eta \,d\mm, \quad &\text{if }\mu\ll\mm\text{ with }\eta:=\frac{d\mu}{d\mm}\\
+\infty, &\text{else.}
\end{array}\right.$$
%and its upgrade by
%$$\Ent^*(\mu\mid\mm)=\left\{
%\begin{array}{ll}
%-\int_\M \log \rho \,d\mu, \quad &\text{if }\mm\ll\mu\text{ with }\rho:=\frac{d\mm}{d\mu}\\
%+\infty. &\text{else.}
%\end{array}\right.$$
It is also known under the name    \emph{Boltzmann entropy} or  \emph{Shannon entropy} or \emph{Kullback-Leibler divergence}.
\end{definition}

%%\begin{lemma} For each $\mu\in\cP(\M)$,
%%$$\Ent(\mu\mid\mm)\le \Ent^*(\mu\mid\mm),$$
%%with equality if $\Ent^*(\mu\mid\mm)<\infty$.
%%\end{lemma}
%\begin{lemma} For every $\mu\in\cP(\M)$ with $\Ent^*(\mu\mid\mm)<\infty$,
%$$\Ent(\mu\mid\mm)=\Ent^*(\mu\mid\mm).$$
%%with equality if $\Ent^*(\mu\mid\mm)<\infty$.
%\end{lemma}
%\begin{proof} 
%%Let $\Ent(\mu\mid\mm)$ be defined as in Definition \ref{def-ent} and denote by $I(\mu)$ the quantity considered in the Lemma. If  $\Ent(\mu\mid\mm)=\infty$ and $I(\mu)=\infty$, the claim obviously holds.
%
%Assume $\Ent^*(\mu\mid\mm)<\infty$. Then $\mu\ll\mm$. (Otherwise, $\exists A$ s.t.~$\mu(A)>0=\mm(A)$. Hence, $\frac{d\mm}{d\mu}=0$  $\mu$-a.e.~on $A$ and thus $\Ent^*(\mu\mid\mm)=\infty$.)
%Moreover, $0<\frac{d\mm}{d\mu}<\infty$ $\mu$-a.e. Thus
%\begin{align*}
%\Ent^*(\mu\mid\mm)&=\int_\M \log \frac{d\mu}{d\mm} \,d\mu=\int_\M \log \frac{d\mu}{d\mm} \, \frac{d\mu}{d\mm}\,d\mm=\Ent(\mu\mid\mm).
%\end{align*}
%%On the other hand, assume $\Ent(\mu\mid\mm)<\infty$. Then $\mu\ll\mm$ and $0<\frac{d\mu}{d\mm}<\infty$ $\mm$-a.e. Thus
%%\begin{align*}
%%\Ent(\mu\mid\mm)=\int_{\{\frac{d\mu}{d\mm}>0\}} \log \frac{d\mu}{d\mm} \, d\mu=
%%-\int_{\{\frac{d\mu}{d\mm}>0\}} \log \frac{d\mm}{d\mu} \, d\mu=
%%-\int \log \frac{d\mm}{d\mu} \, d\mu=I(\mu).
%%\end{align*}
%%{\color{red} Attention: $ \frac{d\mm}{d\mu} $ might not exist!}
%\end{proof}

A key observation now is that under the conjugation map, the role of the measures $\mu$ and $\mm$ in the relative entropy will be exchanged. (Recall that $\mm\cc=\mm$.)
\begin{theorem}
\label{conj-entr}
For every $\mu\in\cP(\M)$,
 \begin{equation*}
 \Ent(\mu \mid \mm) = \Ent(\mm \mid \mu^\c).
\end{equation*}
\end{theorem}
%In {\cite[Cor.~3.4]{Sturm-multi}}, this was stated under more restrictive assumptions.

\begin{proof} 
(i) Let us first remark the well-known fact that both functions $\mu\mapsto  \Ent(\mu \mid \mm)$ and $\mu\mapsto  \Ent(\mm \mid \mu)$ are lower semicontinuous, and secondly observe that both of them are decreasing under the heat flow $t\mapsto \mu_t:=P_t\mu=\eta_t\,\mm$.
Indeed,
\begin{align*}
 \Ent(\mu_t \mid \mm)- &\Ent(\mu \mid \mm)=-4\int_0^t \int_\M |\nabla \sqrt\eta_t|^2\le 0, \\
   \Ent(\mm\mid \mu_t)= &\Ent(\mm \mid \mu)=-\int_0^t \int_\M |\nabla \log\eta_t|^2\le 0.
   \end{align*}
 Thus in particular
 \begin{equation}\label{ent-limits}
 \lim_{t\to0}\Ent(\mu_t \mid \mm)=\Ent(\mu \mid \mm), \qquad \lim_{t\to0}\Ent(\mm \mid \mu_t)=\Ent(\mm \mid \mu).
 \end{equation}

(ii) Observe that $\mu_t$ for each $t>0$ is absolutely continuous w.r.t.~$\mm$ and the density $\eta_t$ is bounded from above and away form 0. 
We claim that
  \begin{equation}\label{dual1}
 \Ent(\mm \mid \mu_t )= \Ent((\mu_t)^\c\mid\mm).
\end{equation}
Indeed,
\begin{align*} \Ent(\mm \mid \mu_t )= \int \log\frac1{\eta_t}\,d\mm=\int \log{\rho_t\circ f_t}\,d\mm=\int\log\rho_t\, d(\mu_t)^\c=\Ent((\mu_t)^\c\mid\mm)
\end{align*}
(where %for convenience we dropped the $t$-dependence from 
the notation follows
Lemma \ref{densities}, in particular,
$\mu_t=\eta_t\mm$, $(\mu_t)^\c=(f_t)_*\mm=\rho_t\mm$).

Taking into account the second assertion in \eqref{ent-limits}, the lower semicontinuity of $\mu\mapsto  \Ent(\mu \mid \mm)$, and the fact that $(\mu_t)^c\to \mu^c$, we obtain
\begin{align}\label{upper-ent} \Ent(\mm\mid \mu)&=\lim_{t\to0}\Ent(\mm \mid \mu_t)=\lim_{t\to0}\Ent( (\mu_t)^c\mid\mm)\ge \Ent(\mu^\c \mid \mm).
\end{align}

(iii) Our next claim is
  \begin{equation*}
 \Ent(\mu_t \mid \mm) = \Ent(\mm \mid (\mu_t)^\c).
\end{equation*}
This can be proven in the same manner as the previous \eqref{dual1}. Alternatively, we can refer to  {\cite[Cor.~3.4]{Sturm-multi}}.
Taking into account the first assertion in \eqref{ent-limits}, the lower semicontinuity of $\mu\mapsto  \Ent(\mm \mid \mu)$, and the fact that $(\mu_t)^c\to \mu^c$, we now obtain
\begin{align}\label{lower-ent} \Ent(\mu\mid \mm)&=\lim_{t\to0}\Ent(\mu_t \mid \mm)=\lim_{t\to0}\Ent(\mm \mid (\mu_t)^c)\ge \Ent(\mm \mid \mu^c).
\end{align}

(iv) Applying this  estimate to $\mu^\c$ in the place of $\mu$ yields
$$
 \Ent(\mu^\c \mid \mm) \ge\Ent(\mm\mid \mu).
$$
Together with \eqref{upper-ent}, this proves the claim.
\end{proof}

\subsection{Support Properties of Conjugate Measures}

\begin{definition} We say that the Riemannian compactum $\M$ is \emph{conjugation regular} if 
for every $z\in\M, \rho\in\cP(\M)$ and $\lambda\in(0,1)$
 there exists an open set $U\subset\M$ with $\mm(U)=\lambda$ such that the optimal transport map $T$ for the transport from $\mm$ to $$\mu:=\lambda\delta_z+(1-\lambda)\rho$$ satisfies
$T(x)=z$ for every $x\in U$, and the conjugate measure $\mu^\c$ does not charge $U$, i.e. $$\mu^\c(U)=0.$$
\end{definition}

\begin{theorem} The following spaces are conjugation regular: 
\begin{itemize}
\item every convex Euclidean compactum,
\item the $n$-torus,
\item the $n$-sphere.
\end{itemize}
\end{theorem}

\begin{proof} For convex Euclidean compacta, see \cite[Lemma 4.1]{Sturm-multi}.
For the torus and the sphere, the claim will be proven in the next subsections as Propositions \ref{OCP-torus} and \ref{OCP-sphere}.
\end{proof}

\begin{conjecture} Every closed Riemannian  manifold is conjugation regular
\end{conjecture}

\begin{theorem} \label{thm1} Assume that $\M$ is conjugation regular. Then for every $\mu\in\cP(\M)$:
\begin{enumerate}[\rm (i)]
\item If $\mu$  is discrete then the topological support of $\mu^\c$  is a  $\mm$-zero set. In particular, $\mu^\c$ has no  $\mm$-absolutely continuous part.

\item If $\mu$ is has full topological support then $\mu^\c$  has no atoms.
\end{enumerate}
\end{theorem}

\begin{proof} (i) Assume that $\mu$ is discrete, say $\mu=\sum_{i=1}^\infty \lambda_i\delta_{z_i}$.
Then according to the previous Lemma,  there exist disjoint open sets $U_i\subset\M$ with $\mm(U_i)=\lambda_i$
and
the conjugate measure $\mu^\c$ is supported by the closed $\mm$-null set $\M\setminus \bigcup_i U_i$. Thus it can not have an absolutely continuous part.

(ii)
To prove the second assertion, assume that $\mu^\c$ has an atom, say $\mu^\c=\lambda\delta_z+(1-\lambda)\rho$. Then according to the  previous Lemma (applied to $\mu^\c$ in the place of $\mu$),
 an open set $U\subset\M$ with $\mm(U)=\lambda$ exists
such that 
 the conjugate measure $\mu=(\mu^\c)^\c$ is supported on  $ \M\setminus  U$ and thus 
 does  not have full  topological support.
\end{proof}

\subsection{Conjugation Regularity of the Torus}
Throughout this subsection, let $\M=\R^n/\Z^n$ be the $n$-torus, occasionally identified with the set $[-1/2,1/2]^n\subset\R^n$.
\begin{lemma} 
Given two points $z,y\in \R^n/\Z^n$ and $\gamma\in\R$ put
\begin{align*}U&:=\left\{x\in \R^n/\Z^n: \d^2(x,z)< \d^2(x,y)+\gamma\right\}\\
A&:=\left\{x\in \R^n/\Z^n: \d^2(x,z)\le \d^2(x,y)+\gamma\right\}.
\end{align*}
Then there exists a convex polyhedron $K\subset [-1/2,1/2]^n\subset\R^n$ such that
$$\pi(K)=A, \quad \pi(K^o)=U$$
 under the quotient map $\pi: \R^n\to\R^n/\Z^n$ with $\pi(0)=z$.
\end{lemma}

Note that the set $A$ is not necessarily convex w.r.t.~the metric on the torus.

\begin{proof} Without restriction $z=0$ and $y\in [0,\frac12]^2$. Indeed, we even assume $y\in (0,\frac12]^2$ and leave the remaining case as an exercise to the reader.

The projection map will be suppressed in the notation, that is, we identify points in the torus with points in $\R^n$. To express the torus-distance in terms of the Euclidean distance, we need $2^n$ copies of $y$:
$$y^\sigma:=y-\sigma, \qquad \sigma\in\{0,1\}^n,$$
and we decompose $[-1/2,1/2]^n$ into $2^n$  subsets (disjoint up to $\mm$-zero sets):
$$R_\sigma:=\prod_{i=1}^n J_{i,\sigma_i}, \qquad J_{i,0}:=[-1/2,-1/2+y_i], \quad J_{i,1}:=[-1/2+y_i,1/2].$$
Then
$$\d(x,y)=
|x-y^\sigma| \qquad \forall x\in R_\sigma,  \ \sigma\in\{0,1\}^n$$
whereas
$$\d(x,z)=|x-z|\qquad \forall x\in [-1/2,1/2]^n.$$
Furthermore, define  closed half-spaces
$$H_\sigma:=\left\{x\in\R^n: \ 
2 x\cdot y^\sigma \le |y^\sigma|^2+\gamma\right\}$$
with boundaries given by $(n-1)$-dimensional hyperplanes
$$L_\sigma:=\left\{x\in\R^n: \ 
2 x\cdot y^\sigma=|y^\sigma|^2+\gamma\right\}.$$
Then for $\sigma\in\{0,1\}^n$ and $x\in R_\sigma$,
\begin{align*}
x\in K\ \Longleftrightarrow \ x\in H_\sigma,\qquad\quad
x\in \partial K\ \Longleftrightarrow \ x\in L_\sigma.
\end{align*}
Thus the set $K$ is the union of the $2^n$ convex polyhedrons
$K_\sigma:=R_\sigma\cap H_\sigma$, $\sigma\in\{0,1\}^n$.

\smallskip

It remains to prove that the set $K$ is convex. The crucial argument for doing so, is that the for any pair of `neighboring' polyhedrons 
$R_\sigma$ and $R_\tau$ with separating hyperplane $V_i$, the hyperplane $L_\sigma$ intersects $V_i$ in exactly the same $(n-2)$-plane as does the hyperplane $L_\tau$. More in detail:
For $i=1,\ldots,n$ define  $(n-1)$-dimensional hyperplanes by
$V_i:=\{x: x_i=y_i-1/2\}$. 
Then obviously $R_\sigma\cap R_\tau=V_i$
for all $\sigma,\tau\in\{0,1\}^n$ with $\sigma_i\not=\tau_i$ and $\sigma_j=\tau_j$ ($\forall j\not=i$),
and straightforward calculations show that
\begin{itemize}
\item[$\bullet$] $L_\sigma\cap L_\tau\subset V_i$
\item[$\bullet$] $L_\sigma\cap V_i=L_\tau\cap V_i=L_\sigma\cap L_\tau.$
\end{itemize}
Thus the set $K_\sigma\cap K_\tau$ is convex.
Iterating this argument to deal with chains of neighboring $K_\sigma$'s (and making use of a continuity argument to deal with transitions between non-neighboring $K_\sigma$'s)  yields the convexity of $K_\sigma\cap K_\tau$ for arbitrary $\sigma,\tau\in\{0,1\}^n$ and thus finally the convexity of 
$$K=\bigcup_{\sigma\in\{0,1\}^n} K_\sigma.$$
\end{proof}

\begin{proposition} \label{OCP-torus}
Let $\M$ be the $n$-torus  and let
$\mu=\lambda\delta_z+(1-\lambda)\rho$ for some $z\in\M, \rho\in\cP(\M)$ and $\lambda\in(0,1)$. 

 Then there exists an open set $U\subset\M$ with $\mm(U)=\lambda$ such that the optimal transport map $T=\exp(\nabla\varphi)$ for the transport from $\mm$ to $\mu$ is well-defined and unique at every point of $U$ and satisfies
$T(x)=z$ for every $x\in U$, and the conjugate measure $\mu^\c$ does not charge $U$, i.e. $$\mu^\c(U)=0.$$
\end{proposition}

\begin{proof} 
(i) Without restriction, we may assume that $\rho(\{z\})=0$. Let $\varphi$ be the optimal $\c$-convex potential (unique up to additive constants) for the transport from $\mm$ to $\mu$ and let
$\Gamma:=\big\{(x,y): \ -\varphi(x)-\varphi^\c(y)=\frac12 \d^2(x,y)\big\}$
denote the contact set (or $\c$-subdifferential).
Let $Y$ be a countable dense subset of $\M\setminus\{z\}$.
Then by continuity
$$\varphi(x)=\sup_{y\in \M}\varphi_y(x)=\sup_{y\in Y}\varphi_y(x)$$
where
$\varphi_y(x):=-\varphi^\c(z)-\frac12\d^2(x,y)$.

Then a necessary condition for a point $x\in\M$ to be transported to $z$ is  lying in
$$K:=\big\{ \varphi_z=\varphi\big\}=
\bigcap_{y\in Y}\{\varphi_z\ge\varphi_y\}= \bigcap_{y\in Y} K_y.
$$
 On the other hand, 
 a necessary condition for a point $x\in\M$ to be transported to $z$ is  lying in
\begin{align*}
\tilde U:= \bigcap_{y\in Y}\{\varphi_z>\varphi_y\}&= \bigcap_{y\in Y} U_y.
\end{align*}
According to the previous Lemma
\begin{align*}
U_y:=&\{\varphi_z>\varphi_y\}=\{\d^2(x,z)<\d^2(x,y)+2\varphi^\c(y)-2\varphi^\c(z)\}\\
K_y:=&\{\varphi_z\ge\varphi_y\}=\{\d^2(x,z)\le \d^2(x,y)+2\varphi^\c(y)-2\varphi^\c(z)\}
\end{align*}
are (projections of) convex ployhedrons in $\R^n$ (open or closed, resp.).
Since $\mm(K_y\setminus U_y)=0$ for all $y\in Y$,
we have $\mm(K\setminus \tilde U)=0$ and thus
$$\mm(K)=\lambda.$$
Since $K$ (regarded as subset of $\R^n$) is convex, $\mm(\partial K)=0$ and thus $U:=K\setminus \partial K$ is the requested convex open set with $\mm(U)=\lambda$ and such that
the optimal transport map $T$ for the transport from $\mm$ to $\mu$  maps $U$ onto $z$. 
More precisely,  we firstly conclude that $T(x) =\exp_x(\nabla\varphi(x))=z$ 
for $\mm$-a.e.~$x\in U$. Since $\varphi$ is semiconvex, it then follows that $T(x) =z$ for all $x\in U$.

(ii) By general principles, the contact set $\Gamma$ is a $\c$-cyclically monotone set. Now assume that $\exists y\not=z,  x\in U$ such that $(x,y)\in\Gamma$. In order to derive a contradiction from this, we again use the geometric concepts and notations from the previous Lemma.
Without restriction, we assume  $z=0$, $y\in (0,\frac12]^2$,
and decompose $[-1/2,1/2]^n$ into the $2^n$  subsets
$R_\sigma:=\prod_{i=1}^n J_{i,\sigma_i}$  (disjoint up to $\mm$-zero sets) such that
$\d(x,y)=
|x-y^\sigma|$ for $x\in R_\sigma$.

Given $z,y,\sigma$ and $x$ as above, 
choose $w\in U\cap R_\sigma$ such that $\langle x-w,z-y\rangle>0$ (which is always possible since $U$ is open). Then 
\begin{align*}\Big[\d^2(x,y)+\d^2(w,z)\Big]-\Big[\d^2(x,z)+\d^2(w,y)\Big]=\langle x-w,y-z\rangle<0\end{align*}
which contradicts the $\c$-cyclic monotonicity and thus the optimality of the coupling.
 Thus for all $x\in U$,
$$\{y: (x,y)\in \Gamma\}=\{z\}.$$
Since the optimal coupling $q$ for the transport from $\mu^\c$ to $\mm$ is supported by $\Gamma$, it follows
$$q\big(U\times (\M\setminus\{z\})\big)\le q\big( \M^2\setminus \Gamma\big)  =0.$$
Furthermore, 
$q\big(U\times \{z\})\big)\le q(\M\times\{z\})=\mm(\{z\})=0$ and thus finally
$$\mu^\c(U)=q(U\times\M)=0.$$
\end{proof}

\subsection{Conjugation Regularity of the Sphere}
Throughout this subsection, let $\M=\mathbb S^n$ be the $n$-torus  identified with the set $\{x\in \R^{n+1}: |x|=1\}$, fix $z\in\M$, and denote its antipodal point by $z^\c$. Put $\M_\epsilon:=B_{\pi-\epsilon}(z)=\M\setminus \overline B_{\epsilon}(z^\c)$.

\begin{lemma} For every $\epsilon >0$ there exists a constant $C=C_\epsilon$ such that for all $y\in \mathbb S^n$ the function $V(x):=\d^2(x,y)-\d^2(x,z)$ satisfies 
\begin{itemize}
\item[\rm(i)]
$\frac1C\, \d(y,z)\le|\nabla_x V(x)|\le 4\pi$ for $x\in \M_\epsilon.$
\item[\rm(ii)]
$\langle\xi,\nabla_{xx} V(x) \xi\rangle\ge -C\,\d(y,z)$ for $x\in \M_\epsilon$ and $\xi\in T_x\mathbb S^n$ with $|\xi|=1$.
\end{itemize}
\end{lemma}

\begin{proof}
(i) The upper bound follows trivially from $\d(x,.)|\le\pi$. For the lower bound, we consider three cases:

(a) If $x\in B_{\pi-\epsilon}(z)\cap B_{\pi-\epsilon/4}(y)$ then the lower bound with $c_\epsilon:= \frac{\pi-\epsilon/4}{\sin(\pi-\epsilon/4)}$ (for $\epsilon<\pi/2$) follows from distance comparison for the two geodesics emanating from $x$ with tangents $\nabla_x \d^2(x,y)$ and $\nabla_x \d^2(x,z)$, resp.

(b) If $\d(y,z)\le \frac34\epsilon$ (or equivalently, $\d(y^\c,z^\c)\le\frac34\epsilon$) and $\d(x,z^\c)>\epsilon$ then $\d(x,y^\c)>\frac14\epsilon$  and thus we are back in case (a).

(c) If If $\d(y,z)>\frac34\epsilon$ (or equivalently, $\d(y^\c,z^\c)>\frac34\epsilon$) and $\d(x,y^\c)\le\frac14\epsilon$ then $\d(x,z^\c)>\frac12\epsilon$.  Therefore
\begin{align*}
|\nabla_x\d^2(x,y)|&=2\,|\d(x,y)|\ge 2(\pi-\frac14\epsilon), \\
|\nabla_x\d^2(x,z)|&=2\,|\d(x,z)|< 2(\pi-\frac12\epsilon)
\end{align*}
and thus
$$|\nabla_x\d^2(x,y)-\nabla_x\d^2(x,z)|> \frac14\epsilon.$$

(ii) Recall that the Hessian $\nabla_{xx} \d^2(x,w)$ of the squared distance function has eigenvalues 2 and $2\frac{r\cos r}{\sin r}\in (-\infty,2]$ with $r:=\d(x,w)$, the latter with multiplicity $n-1$.
Thus for every $y\in\M$, every $x\in\M_\epsilon$,  and every  $\xi\in T_x\mathbb S^n$ with $|\xi|=1$,
\begin{align*}\langle\xi,(\nabla_{xx} \d^2(x,z)- \nabla_{xx}\d^2(x,y))\, \xi\rangle\ge 2\frac{\d(x,z)\,\cos\d(x,z)}{\sin\d(x,z)}-2\ge -2\pi \, \cot\epsilon=:-C^*.
\end{align*}
Moreover,  for every $x\in\mathbb S^n\setminus\{z^\c\}$ and every $\xi\in T_x\mathbb S^n$,
$$\langle\xi,\nabla_{xx} \d^2(x,z)\, \xi\rangle=\langle\xi^z,\nabla_{zz} \d^2(x,z)\, \xi^z\rangle$$
where $\xi^z\in T_z\mathbb S^n$ denotes the parallel transport of $-\xi$ along the minimal geodesic connecting $x$ and $z$. Analogously, with $y$ in the place of $z$.
Thus for 
 $x\in\mathbb S^n\setminus\{z^\c, y^\c\}$
\begin{align*}\langle\xi,(\nabla_{xx} \d^2(x,z)&- \nabla_{xx}\d^2(x,y))\, \xi\rangle\\
&=\langle\xi^z,\nabla_{zz} \d^2(x,z)\, \xi^z\rangle-
\langle\xi^y, \nabla_{yy}\d^2(x,y)\, \xi^y\rangle.
\end{align*}
For $y$ close to $z$ we can parallel transport the vector $\xi^y$ and the bilinear form  $\nabla_{yy}$ to $T_z\mathbb S^n$
and thus obtain
\begin{align*}\langle\xi,(\nabla_{xx} \d^2(x,z)&- \nabla_{xx}\d^2(x,y))\, \xi\rangle\\
&=\langle\xi^z-\tilde\xi^y,\nabla_{zz} \d^2(x,z)\xi^z+\tilde\nabla_{yy}\d^2(x,y) \tilde\xi^y\rangle\\
&\quad+
\langle\xi^z, (\nabla_{zz} \d^2(x,z)-\tilde\nabla_{yy}\d^2(x,y))\, \tilde\xi^y\rangle.
\end{align*}
For  $y\in B_{\epsilon/2}(z)$ and  $x\in \M_\epsilon$ (which  implies $\d(x,y^\c)>\epsilon/2$), the eigenvalues of $\nabla_{zz} \d^2(x,z)$ and $\tilde\nabla_{yy}\d^2(x,y)$ are bounded in modulus by $2\pi \, \cot(\epsilon/2)$ and the eigenvalues of $\nabla_{zz} \d^2(x,z)-\tilde\nabla_{yy}\d^2(x,y)$ are bounded in modulus by $C \,\d(x,y)$ for some constant $C$.
Furthermore, $|\xi^z-\tilde\xi^y|\le C' \,\d(x,y)$. Thus
\begin{align*}\langle\xi,(\nabla_{xx} \d^2(x,z)&- \nabla_{xx}\d^2(x,y))\, \xi\rangle
\ge - C'' \, \d(x,y)\end{align*}
for every $\xi\in T_x\mathbb S^n$ with $|\xi|=1$.
\end{proof}

\begin{lemma} 
For every $\epsilon >0$ there exists a constant $C=C_\epsilon$ such that for every $y\in \mathbb S^n$ and every $\gamma\in\R$ the boundary $\partial U$ of the set
\begin{align*}U&:=U_\epsilon(y):=\left\{x\in {\mathbb S}^n: \d^2(x,z)< \d^2(x,y)+\gamma\right\}\ \cap \ \M_\epsilon
\end{align*}
has  second fundamental form  bounded from below by $-C$:
$$\Pi_{\partial U}(\xi,\xi)\ge -C \, |\xi|^2 \qquad \forall \xi \in T\partial U.$$
\end{lemma}

\begin{proof} Given $z$ and $y$ as above, put $V(x)= \d^2(x,y)- \d^2(x,z)+\gamma$. Then 
$$\partial U=\left(\{V=0\}\cap \M_\epsilon\right) \cup \left(\{V\ge0\}\cap\partial \M_\epsilon\right).$$
For $x\in \{V=0\}\cap \M_\epsilon$ and $\xi \in T_x\partial U$ with $|\xi|=1$,
$$\Pi_{\partial U}(\xi,\xi)=\frac{\langle\xi,\nabla^2 V(x) \xi\rangle}{|\nabla V(x)|}\ge -C^2$$
according to the previous Lemma.
Furthermore, the second  fundamental form of $\partial \M_\epsilon=\partial B_\epsilon(z^\c)$ is bounded from below by $-\frac1\epsilon$.
Thus the claim follows by the $\cap$-stability of lower bounds for the second  fundamental form of boundaries.
\end{proof}

\begin{proposition} \label{OCP-sphere}
Let $\M$ be the $n$-sphere  and let
$\mu=\lambda\delta_z+(1-\lambda)\rho$ for some $z\in\M, \rho\in\cP(\M)$ and $\lambda\in(0,1)$. 

 Then there exists an open set $U\subset\M$ with $\mm(U)=\lambda$ such that the optimal transport map $T=\exp(\nabla\varphi)$ for the transport from $\mm$ to $\mu$ is well-defined and unique at every point of $U$ and satisfies
$T(x)=z$ for every $x\in U$, and the conjugate measure $\mu^\c$ does not charge $U$, i.e. $$\mu^\c(U)=0.$$
\end{proposition}

\begin{proof} We follow the argumentation of the proof for Proposition \ref{OCP-torus}. 

(i) Again,
 a necessary condition for a point $x\in\M$ to be transported to $z$ is  lying in
$$K:=\big\{ \varphi_z=\varphi\big\}=
\bigcap_{y\in Y}\{\varphi_z\ge\varphi_y\}= \bigcap_{y\in Y} K(y),
$$
and 
 a necessary condition for a point $x\in\M$ to be transported to $z$ is  lying in
\begin{align*}
\tilde U:= \bigcap_{y\in Y}\{\varphi_z>\varphi_y\}&= \bigcap_{y\in Y} U(y)
\end{align*}
where
\begin{align*}
U(y):=&\{x\in\mathbb S^n: \ \d^2(x,z)<\d^2(x,y)+2\varphi^\c(y)-2\varphi^\c(z)\}\\
K(y):=&\{x\in\mathbb S^n: \ \d^2(x,z)\le \d^2(x,y)+2\varphi^\c(y)-2\varphi^\c(z)\}.
\end{align*}
According to the previous Lemma, for every $\epsilon>0$ each of the sets
$$U_\epsilon(y):=U(y)\cap \M_\epsilon$$
is semiconvex: the second fundamental form of $\partial U_\epsilon(y)$ is bounded from below by $-C_\epsilon$, uniformly in $y\in\M$.
Thus also $\tilde U_\epsilon:=\tilde U\cap \M_\epsilon$ is semiconvex with  second fundamental form of its boundary bounded from below by $-C_\epsilon$.  In particular, $\mm(\partial\tilde U_\epsilon)=0$. Let $U_\epsilon$ denote the interior of  $\tilde U_\epsilon$. Then
$$U:=\bigcup_{\epsilon>0}U_\epsilon$$
is the requested open set with $\mm(U)=\lambda$ and the property that every point $x\in U$ is transported to $z$.

(ii)
Given $z,y$ and $x$ as above, 
assume without restriction that $z$ and $y$ lie in the $e_1$-$e_2$-plane, say,
$y=(1,0,\ldots, 0)$, $z=(\cos\varphi, \sin\varphi, 0,\ldots,0)$
with $\varphi\in(0,\pi)$. Then for sufficiently small $t>0$ the point $w:=\Phi_{t}(x)$ lies in $U$ where
$\Phi_t: \big(r \cos\psi, r\sin\psi, y_3, \ldots,y_n\big) \mapsto \big(r \cos(\psi+t\varphi), r\sin(\psi+ t \varphi), y_3, \ldots,y_n\big)$
denotes rotation of the $e_1$-$e_2$-plane.
Then 
\begin{align*}\Big[\d^2(x,y)+\d^2(w,z)\Big]-\Big[\d^2(x,z)+\d^2(w,y)\Big]<0\end{align*}
which contradicts the $\c$-cyclic monotonicity and thus the optimality of the coupling.
 Therefore, for all $x\in U$,
$$\{y: (x,y)\in \Gamma\}=\{z\}$$
and thus finally
$\mu^\c(U)=q(U\times\M)=0$.
\end{proof}

\section{Entropic Measure}
\subsection{Heuristics}

\paragraph{(a)}
Our goal is to construct a canonical probability measure $\mathbb{P}^\beta$ on the McCann--Otto  space
$\mathcal{P} =\mathcal{P}(\M)$ over a closed Riemannian manifold, according to the formal ansatz
\begin{equation}\label{heur}
\mathbb{P}^\beta(d\mu) =\frac{1}{Z} e^{-\beta \, \Ent(\mu \mid \mm)}\ \mathbb{P}^*(d\mu).
\end{equation}
Here $\Ent(\cdot \mid \mm)$ is the \emph{relative entropy} with respect to the normalized volume measure $\mm$ and $\beta$ is a constant $>0$ (`the inverse temperature');
$\mathbb{P}^*$ should denote a (non-existing) `uniform distribution' on $\mathcal{P}(\M)$ and $Z$  a normalizing constant. Using the conjugation map
$\Conj_\mathcal{P}: \mathcal{P}(\M)\rightarrow \mathcal{P}(\M)$ --- which will be introduced and analyzed in the next subsection --- and denoting $\mathbb{Q}^\beta := (\Conj_\mathcal{P})_\ast \mathbb{P}^\beta$, $\mathbb{Q}^* :=
(\Conj_\mathcal{P})_\ast \mathbb{P}^*$, the above problem can be reformulated as follows:

Construct a probability measure $\mathbb{Q}^\beta$ on $\mathcal{P}(\M)$ such that --- at least formally ---
\begin{equation}{\label{Qbeta}}
 \mathbb{Q}^\beta (d\nu) = \frac{1}{Z} e^{-\beta \, \Ent (\mm \mid \nu)} \  \mathbb{Q}^*(d\nu)
\end{equation}
with some `uniform distribution' $\mathbb{Q}^*$ in $\mathcal{P}(\M)$.
Here, we have used the fact that $$\Ent (\nu\cc \mid \mm) = \Ent (\mm \mid \nu)$$ (Prop.~\ref{conj-entr}), at least if $\nu \ll \mm$ with $\frac{d \nu}{d\mm}>0$ almost everywhere.

\paragraph{(b)}

 The representation (\ref{Qbeta}) is reminiscent
of Feynman's heuristic picture of the Wiener measure. Let us briefly
recall the latter and try to use it as a guideline for our
construction of the measure $\Q^\beta$.

According to this heuristic picture, the Wiener measure
$\mathbf{P}^\beta$ with diffusion constant $\sigma^2=1/\beta$ should
be interpreted (and could be constructed) as
\begin{equation}\label{Feynman}
\mathbf{P}^\beta(dg)=\frac1{Z_\beta}\, e^{-\beta\cdot H(g)}\,
\mathbf{P}^*(dg)
\end{equation} with the energy functional $H(g)=\frac12\int_0^1
g'(t)^2dt$. Here $\mathbf{P}^*(dg)$ is assumed to be the 'uniform
distribution' on the space $\C_0=\C_0([0,1],\R)$ of all continuous paths
$g: [0,1]\to\R$ with $g(0)=0$. Even if such a uniform distribution existed, typically almost all paths $g$ would have infinite energy. Nevertheless, one can overcome this difficulty as follows.

Given any finite partition $\{0=t_0<t_1< \dots< t_N=1\}$ of $[0,1]$,
one should replace the energy $H(g)$ of the path $g$ by $H_N(g)=\inf\big\{ H(\tilde g): \ \tilde g\in\C_0, \ \tilde g(t_i)=g(t_i) \
\forall i\big\}$ which coincides with the energy
of the piecewise linear interpolation of $g$,
$$H_N(g)=
\sum_{i=1}^N\frac{|g(t_i)-g(t_{i-1})|^2}{2(t_i-t_{i-1})} .
$$
Then (\ref{Feynman}) leads to the following explicit representation
for the finite dimensional distributions
\begin{equation}\label{Feynman2}
\mathbf{P}^\beta\left(g_{t_1}\in dx_1,\ldots,g_{t_N}\in dx_N\right)=
\frac1{Z_{\beta,N}}\exp\left(-\frac\beta2
\sum_{i=1}^N\frac{|x_i-x_{i-1}|^2}{t_i-t_{i-1}} \right)\, p_N(dx_1,
\ldots,x_N).
\end{equation}
Here $p_N(dx_1, \ldots,x_N)=\mathbf{P}\left(g_{t_1}\in
dx_1,\ldots,g_{t_N}\in dx_N\right)$ should be a 'uniform
distribution' on $\R^N$ and $Z_{\beta,N}$ a normalization constant.
Choosing $p_N$ to be the $N$-dimensional Lebesgue measure makes the
RHS of (\ref{Feynman2})  a projective family of probability
measures. According to Kolmogorov's extension theorem this family
has a unique projective limit, the Wiener measure $\mathbf{P}^\beta$
on $\C_0$ with diffusion constant $\sigma^2=1/\beta$.

\paragraph{(c)}
Probability measures ${\bf{P}} (d\mu)$ on $\mathcal{P}(\M)$ --- so called \emph{random probability measures} on $\M$ --- are uniquely determined by the distributions
${\bf{P}}_{M_1,\ldots,M_N}$ of the random vectors $$(\mu(M_1),\ldots,\mu(M_N))$$ for all $N\in \mathbb{N}$ and all measurable partitions 
of $\M$ into disjoint measurable subsets $M_i$.
Conversely, if a consistent family ${\bf{P}}_{M_1,\ldots,M_N}$ of probability measures on $[0,1]^N$ (for all $ N\in \mathbb{N}$ and all measurable partitions ${\dot\cup}_{i=1}^N M_i =\M$) is given then there exists a random probability measure ${\bf P}$ such that $$ {\bf P}_{M_1,\ldots,M_N} (A)= {\bf{P}}((\mu(M_1),\ldots,\mu(M_N)) \in
A)$$ for all measurable $A \subset [0,1]^N$, all $ N\in \mathbb{N}$ and all partitions $\dot{\cup}_{i=1}^N M_i=\M$.

Given a measurable partition $\dot{\cup}_{i=1}^N M_i=\M$, the heuristic formula (\ref{Qbeta}) yields the following ansatz for the finite dimensional
distribution on $[0,1]^N$:
\begin{equation}{\label{QbetaM}}
 \mathbb{Q}^\beta_{M_1,\ldots,M_N} (d x) = \frac{1}{Z_N} e^{-\beta S_{M_1,\ldots,M_N}(x)} q_{M_1,\ldots,M_N} (dx)
\end{equation}
where
$q_{M_1,\ldots,M_N}(dx) = \mathbb{Q}^* ((\nu(M_1), \ldots , \nu(M_N)) \in dx)$ denotes the distribution of the random vector
$(\nu(M_1), \ldots, \nu(M_N))$ in the simplex 
$\sum_N = \big\lbrace x \in [0,1]^N : \sum_{i=1}^N x_i =1 \big\rbrace $
and $S_{M_1,\ldots,M_N}(x)$ denotes the 
the {minimum} of $\nu\mapsto\Ent(\mm \mid \nu)$ under the constraint
$\nu(M_1)=x_1, \ldots, \nu(M_N)=x_N$, that is,
\begin{equation}
{S}_{M_1,\ldots,M_N}(x) = -\sum_{i=1}^{N} \log \frac{x_i}{\mm(M_i)} \cdot \mm(M_i).
\end{equation}
%
%
%conditional expectation (with respect to $\mathbb{Q}^*$) of $S(\cdot) = \Ent (\mm \mid \cdot \,)$ under the condition
%$\nu(M_1)=x_1, \ldots, \nu(M_N)=x_N$.
%Moreover,  

\paragraph{(d)}
What is the     `uniform distribution' $q_N$ in the simplex $\sum_N$? 
A natural requirement
will be the invariance under merging/subdividing
\begin{align}\label{q-inv}
q_N(dx_1,\ldots, dx_N)
=\Xi^i_*\Big[
q_{N-1}(dx_1,\ldots, dx_{i-1},d\xi,dx_{i+2}, \ldots, dx_N)\otimes
 q_{2}(dy_1,dy_2)\Big]
\end{align}  
where $$\Xi^i: (x_1,\ldots, x_{i-1},\xi,x_{i+2}, \ldots, x_N), (y_1,y_2)\mapsto (x_1,\ldots, x_{i-1},\xi \cdot y_1, \xi \cdot y_2, x_{i+2}, \ldots, x_N).$$
If the $q_N$, $N\in\N$, were probability measures then the
invariance property  admits the following interpretation: 
\begin{itemize}
\item
under
$q_N$, for each $i$ the distribution of the $(N-1)$-tuple
$(x_1,\ldots, x_{i-1},x_i+x_{i+1},x_{i+2}, \ldots, x_N)$ is given by $q_{N-1}$;
\item and under $q_N$, the distribution of the pair $\big(\frac{x_i}{x_i+x_{i+1}}, \frac{x_{i+1}}{x_i+x_{i+1}}\big)$
 is given by $q_2$.
 \end{itemize}
% 
%\begin{lemma}[{\cite{RS09}, Lemma 3.1}]
%There exists no family of {\em probability} measures $(q_N)_{N\in\N}$ with
% property (\ref{q-inv}).
  According to {\cite[Lemma 3.1]{RS09},  the unique family of measures $(q_N)_{N\ge2}$ with continuous densities 
which satisfies property (\ref{q-inv}) is given by
\begin{equation}{\label{qn}}
 q_N(dx)= C^N\cdot \frac{dx_1 \ldots dx_{N-1}}{x_1 \cdot x_2 \cdot \ldots \cdot x_{N-1} \cdot x_N} \cdot \delta_{(1-\sum_ {i=1}^{N-1} x_i)} (d
 x_N).
\end{equation}
for some constant $C\in\R_+$.
%\end{lemma}

\paragraph{(e)}
%
%To get hands on $S_{M_1,\ldots,M_N}(x)$, the \emph{conditional expectation} of $S(\cdot) =\Ent(\mm \mid \cdot \,)$ under the constraint $\nu(M_1)=x_1,
%\ldots, \nu(M_N)=x_N$, we  replace it by 
%\begin{equation}
%\underline{S}_{M_1,\ldots,M_N}(x) = -\sum_{i=1}^{N} \log \frac{x_i}{\mm(M_i)} \cdot \mm(M_i),
%\end{equation}
%the \emph{minimum} of $\nu\mapsto\Ent(\mm \mid \nu)$ under the constraint
%$\nu(M_1)=x_1, \ldots, \nu(M_N)=x_N$.
%Then
Thus
\begin{align*}
  \mathbb{Q}^\beta_{M_1,\ldots,M_N} (dx) &= c \cdot e^{-\beta {S}_{M_1,\ldots,M_N}(x)} q_N(dx) \\
&= \frac{\Gamma(\beta)}{\overset{N}{\underset{i=1}{\prod}} \Gamma (\beta \mm (M_i))} \cdot 
\prod_{i=1}^N x_i^{\beta \cdot \mm (M_i)-1}  \
 \delta_{(1-{\overset{N-1}{\underset{i=1}{\sum}}} x_i)}(dx_N)\, dx_{N-1} \ldots dx_1.
\end{align*}
This, indeed, defines a projective family. Hence, the random probability measure $\mathbb{Q}^\beta$ exists and is uniquely defined. It is the well-known
\emph{Dirichlet-Ferguson measure}. Therefore, in turn, also the random probability measure $\mathbb{P}^\beta = (\Conj_\mathcal{P})_\ast \mathbb{Q}^\beta$ exists
uniquely.

\subsection{Rigorous Definition}

\begin{definition}
 Given any  Riemannian compactum $\M$ with reference measure $\mm$ and any parameter $\beta > 0$, the entropic measure
$$ \mathbb{P}^\beta := (\Conj_\mathcal{P})_\ast \mathbb{Q}^\beta$$
is the push forward of the Dirichlet-Ferguson measure $\mathbb{Q}^\beta$ (with reference measure $\beta \mm$) under the conjugation map
$\Conj_\mathcal{P} :\mathcal{P}(\M) \rightarrow \mathcal{P}(\M)$.
\end{definition}

$\mathbb{P}^\beta$ as well as $\mathbb{Q}^\beta$ are probability measures on the compact space $\mathcal{P}= \mathcal{P}(\M)$ of probability measures on $\M$.

\begin{center}
\includegraphics[width=7cm]{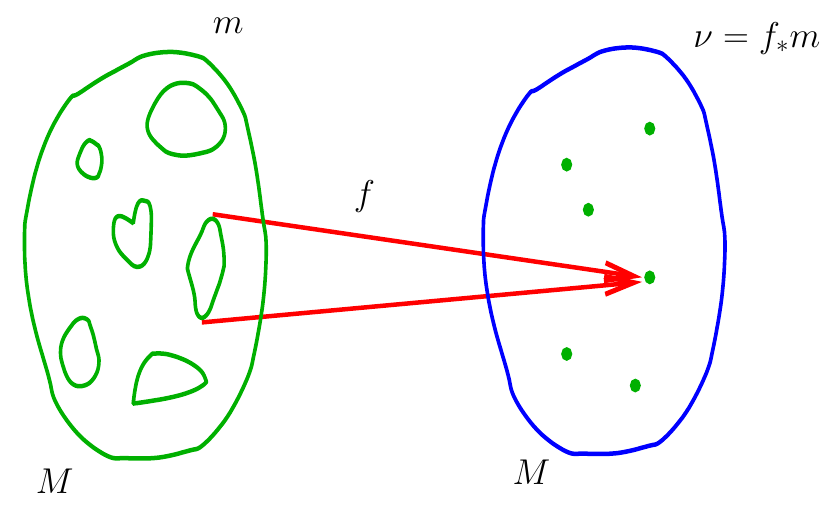}\qquad
\includegraphics[width=7cm]{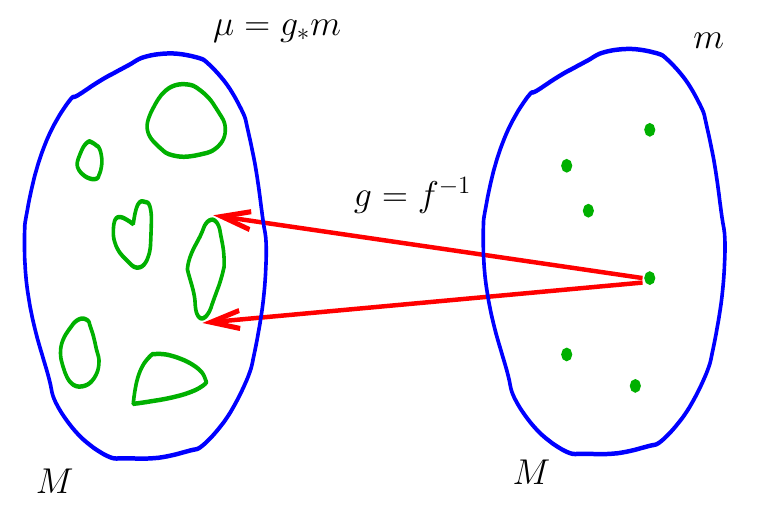}
\end{center}

Recall the definition of the Dirichlet-Ferguson measure $\mathbb{Q}^\beta$ \cite{Ferguson73}: For each measurable
partition $ \dot{\cup}_{i=1}^N M_i=\M$,
the random vector $(\nu(M_1),\ldots,\nu(M_N))$ is distributed according to a Dirichlet distribution with parameters $(\beta \,\mm(M_1),\ldots,\beta\, \mm(M_N))$.
That is, for any bounded Borel function $u: \R^N \rightarrow \R$,
\begin{align}
\label{Ferguson}\int_{\mathcal{P}(M)} u (\nu(M_1),\ldots,\nu(M_N))&\mathbb{Q}^\beta (d\nu) =\\
=\dfrac{\Gamma(\beta)}{{\overset{N}{\underset{i=1}{\prod}}} \Gamma (\beta \mm(M_i))} \cdot \int_{\left[ 0,1\right]^N  } & u(x_1,\ldots,x_N) \cdot x_1^{\beta \, \mm
(M_1)-1} \cdot \ldots \cdot x_N^{\beta \, \mm (M_N)-1} \times\nonumber \\
& \times\delta_{(1-{\overset{N-1}{\underset{i=1}{\sum}}}x_i)}(d x_N) d x_{N-1}\ldots d x_1.\nonumber
\end{align}
The latter uniquely characterizes the `random probability measure' $\mathbb{Q}^\beta$.
The existence (as a projective limit) is guaranteed by Kolmogorov's theorem.

Note that 
$${\mathbb E}_{\Q^\beta}\big[ \nu(A)\big]=\mm(A)$$
for every measurable $A\subset \M$. In particular 
$\mm(A)=0$ implies that $\nu(A)=0$ for $\mathbb{Q}^\beta$-a.e.~$\nu \in \mathcal{P}(\M)$. 
On the other hand, 
$\mathbb{Q}^\beta$-a.e.~$\nu \in \mathcal{P}(M)$ is discrete.
This is obvious by means of an alternative, more direct construction:

Let $(x_i)_{i\in\mathbb{N}}$ be an iid sequence of points in $\M$, distributed according to $\mm$, and let $(t_i)_{i\in\mathbb{N}}$ be an iid sequence of numbers in
$[0,1]$, independent of the previous sequence and distributed according to the Beta distribution with parameters 1 and $\beta$, i.e. ${\textmd{Prob}}(t_i \in ds)=
\beta (1-s)^{\beta-1} \cdot 1_{[0,1]}(s) ds$.
Put 
\begin{equation}\lambda_k = t_k \cdot \prod_{i=1}^{k-1} (1- t_i) \qquad \text{and} \qquad \nu = \sum_{k=1}^{\infty} \lambda_k \cdot \delta_{x_k}.
\label{particle}
\end{equation}
Then $\nu \in \mathcal{P}(\M)$ is distributed according to $\mathbb{Q}^\beta$ \cite{Seth94}.

The distribution of $\nu$ does not change if one replaces the above `stick-breaking process' $(\lambda_k)_{k\in\mathbb{N}}$ by the `Dirichlet-Poisson process' $(\lambda_{(k)})_{k\in\mathbb{N}}$ obtained from it by ordering the entries of the previous one according to their size: $\lambda_{(1)}\ge\lambda_{(2)}\ge\ldots\ge0$.
Alternatively, the Dirichlet-Poisson process can be regarded as the sequence of jumps of a Gamma process with parameter $\beta$, ordered according to size.

For  further properties of the Dirichlet-Ferguson measure we refer to \cite{Lorenzo}.

\subsection{Support Properties}

Let us fix some notation: we denote by $B_r(x)$ the $\d$-balls in $\M$ and by ${\mathbb B}_r(\mu)$ the balls in ${\mathcal P}(\M)$ w.r.t.~the 
 Kantorovich-Wasserstein 2-distance ${\sf W}_2$  derived from $\d$.

\begin{theorem}\label{full-supp} $\PP^\beta$ has full topological support in $\cP(\M)$, that is,
$$\PP^\beta\big(\mathbb B_\epsilon(\mu)\big)>0\qquad \forall \mu,\epsilon.$$
\end{theorem}

\begin{proof} Cf.~\cite[Lem.~3.26]{Lorenzo0}. Since the conjugation map ${\mathfrak C}:\nu\mapsto \nu^{\mathfrak c}$ is continuous, for every $\mu\in\cP(\M)$ and every $\epsilon>0$ there exists a $\delta>0$ such that ${\mathfrak C}\big(\mathbb B_\delta(\mu^\c)\big)\subset \mathbb B_\epsilon(\mu)$. Thus
$$\PP^\beta\big(\mathbb B_\epsilon(\mu)\big)=\Q^\beta\big({\mathfrak C}^{-1}\big(\mathbb B_\epsilon(\mu)\big)\big)
\ge \Q^\beta\big(\mathbb B_\delta(\mu^\c)\big).$$
The latter probability is positive according to the next lemma.
\end{proof}

\begin{lemma} $\Q^\beta$ has full topological support in $\cP(\M)$.
\end{lemma}

\begin{proof} (1) Given $\nu\in\cP(\M)$ and $\epsilon>0$, there exist $k\in\N$, $\alpha_i\in [0,1]$ and $x_i\in\M$ such that
$${\sf W}_2(\nu,\nu')\le \epsilon, \qquad \nu':=\sum_{i=1}^k \alpha_i\delta_{x_i}.$$
Without restriction, $\alpha_1\ge\alpha_2\ge \ldots\ge \alpha_k$.

(2) We use the fact that the random probability measures distributed according to $\Q^\beta$ can be represented as 
$$\nu=\sum_{i=1}^\infty \lambda_i\delta_{z_i}$$
where $(\lambda_i)_{i\in\N}$ denotes a non-increasing sequence of positive numbers distributed according to the Dirichlet-Poisson process with parameter $\beta$, and where  $(z_i)_{i\in\N}$ denotes an iid sequence of uniformly distributed points in $\M$.

There exists $p_1>0$ such that with probability $\ge p_1$ the Dirichlet-Poisson process $(\lambda_i)_{i\in\N}$ with parameter $\beta$ satisfies
$$\lambda_i\in [\alpha_i-\epsilon^2/k,\alpha_i]\qquad \forall i\in\{1,\ldots,k\},$$
and there exists $p_2>0$ such that with probability $\ge p_2$ the iid sequence of uniformly distributed points $(z_i)_{i\in\N}$ satisfies
$$z_i\in B_\epsilon(x_i)\qquad \forall i\in\{1,\ldots,k\}.$$

(3) By transporting each $\lambda_i\delta_{z_i}$ for  $i\in\{1,\ldots,k\}$ to $\lambda_i\delta_{x_i}$ and ``the remainder to the remainder'', that is, $\sum_{i=k+1}^\infty \lambda_i\delta_{z_i}$ to $\sum_{i=1}^k (\alpha_i-\lambda_i)\delta_{x_i}$ we find a coupling which allows us to estimate
\begin{align*}
{\sf W}_2^2\Bigg(\sum_{i=1}^\infty \lambda_i\delta_{z_i}, \sum_{i=1}^k \alpha_i\delta_{x_i}\Bigg)&\le(1-\epsilon^2) \cdot \epsilon^2+ \epsilon^2\cdot L^2\le \epsilon^2\cdot (1+L)^2.
\end{align*}
Here $L:=\max_{x,y\in\M} \d(x,y)$ denotes the diameter of $\M$.

(4) Therefore,
$$\Q^\beta\Big(\mathbb B_{\epsilon(1+L)}(\nu')\Big)\ge p_1\cdot p_2=:p$$
and thus
$$\Q^\beta\Big(\mathbb B_{\epsilon(2+L)}(\nu)\Big)\ge \Q^\beta\Big(\mathbb B_{\epsilon(1+L)}(\nu')\Big)\ge p>0.$$
\end{proof}

\begin{theorem} Assume that $\M$ is conjugation regular. Then  $\PP^\beta$-a.s. 

\begin{itemize}
\item
$\mu$  has no atoms, 
\item $\mu$ has no absolutely continuous part,\\
even more, the topological support  of $\mu$ is a $\mm$-zero set. 

\end{itemize}
\end{theorem}

\begin{proof} 
Immediate consequence of Theorem \ref{thm1} and  the next Lemma which asserts that for $\PP^\beta$-a.e. $\mu$ the measure $\mu^\c$ is discrete and has full topological support. 
\end{proof}

\begin{lemma} $\Q^\beta$-a.e. $\nu$ is discrete and has full topological support in $\M$. \end{lemma}
\begin{proof} As already mentioned before, $\Q^\beta$-a.e.~$\nu$ can be represented as $\sum_{i=1}^\infty \lambda_i\delta_{z_i}$. This proves the discreteness. The support property follows easily from $\nu(B_\epsilon(y))\ge \lambda_1 m(B_\epsilon(y))>0$ for every $y\in\M$, every $\epsilon>0$, and a.e.~$\lambda_1$.
\end{proof}

For  $\mathbb{P}^\beta$-a.e.~$\mu \in \mathcal{P}(\M)$, there exist a countable number of open sets $U_k\subset \M$ (`holes in the support of $\mu$')
with sizes $\lambda_k=m(U_k)$, $k\in\mathbb{N}$. The measure $\mu$ is supported on  the complement of all these holes $\M\setminus\bigcup_k U_k$, a compact $\mm$-zero set.

The sequence $(\lambda_k)_{k\in\mathbb{N}}$ of sizes of the holes is distributed according to the stick breaking process with parameter $\beta$. In particular,
$\mathbb{E}\lambda_k=\frac1\beta(\frac\beta{1+\beta})^k\le\frac1{1+\beta}$ for  $k\in\mathbb{N}$.
For large $\beta$, the size of the $k$-th hole decays like $\frac1\beta\exp(-k/\beta)$ as $k\to\infty$.

\medskip

Typical realizations of $\mu$ as displayed in Figure 1 remind of the simulations of the early universe in \cite{early}:\vspace{-.8cm}
%\begin{figure}
\begin{center}
\includegraphics[width=4cm]{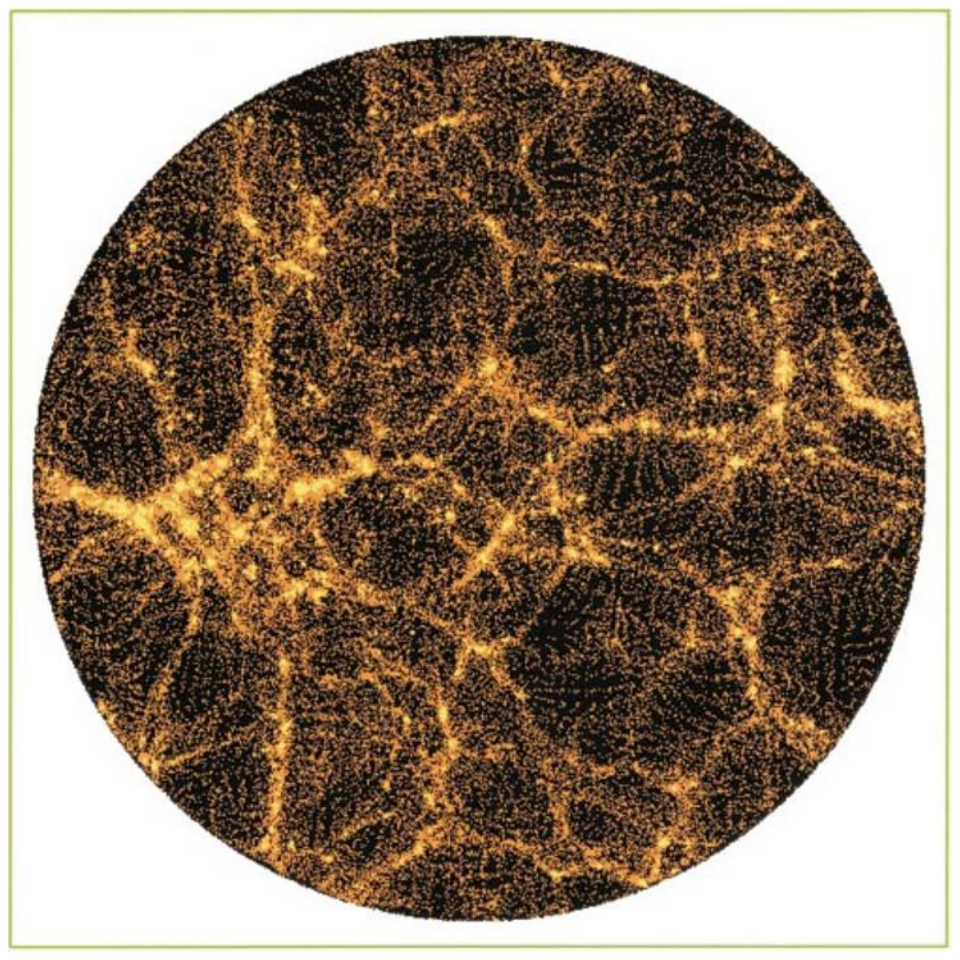}
\end{center}
\vspace{-2.5cm}

\subsection{Mean Value and Asymptotics}
For the sequel, assume that $\M$ is the $n$-sphere or the $n$-torus and that $\mm$ is the respective normalized volume measure.

\begin{theorem} $\PP^\beta$ has mean value $\mm$, i.e.~
$\int_{\cP(\M)}\mu\, d\PP^\beta(\mu)=\mm$ or, more explicitly,
$$\int_{\cP(\M)}\mu(A)\, d\PP^\beta(\mu)=\mm(A)\qquad\quad (\forall \text{Borel } A\subset\M).$$
\end{theorem}

\begin{proof} 
Translation invariance of $\mm$ implies that  the Dirichlet-Ferguson measure is invariant under translations. Furthermore,
translation invariance of $\d$ and $\mm$ implies that  the conjugation map $\mu\mapsto \mu^\c$ commutes with translations. 
Hence, the entropic measure is invariant under translations and so is its mean value. But the only translation invariant probability measure on $\M$ is the given $\mm$.
\end{proof}

\begin{definition}
\begin{enumerate}[\rm (i)]
\item The \emph{uniform distribution on the Dirac measures} is the measure $\PP^0:=\delta_*\mm\in\cP(\M)$ defined as push forward of $\mm$ under the map %(=isometric embedding) 
$\delta: \M\to\cP(\M), x\mapsto \delta_x$. It is characterized by
\begin{align*}\PP^0(A)&=\mm\big\{x\in\M: \delta_x\in A\big\} \qquad  \forall \text{Borel }A\subset \cP(\M),\\
\int_{\cP(\M)} f\,d\PP^0&=\int_\M f(\delta_x)\,d\mm(x) \qquad  \forall \text{Borel } f:\cP(\M)\mapsto\R.
\end{align*}

\item The \emph{Dirac measure on the uniform distribution}  is the measure $\PP^\infty:=\delta_\mm\in\cP(\M)$ characterized by
\begin{align*}\PP^\infty(A)=1_A(\mm)
\qquad \forall \text{Borel }A\subset \cP(\M),\\
 \int_{\cP(\M)} f\,d\PP^\infty= f(\mm)\qquad  \forall \text{Borel } f:\cP(\M)\mapsto\R.
\end{align*}
\end{enumerate}
\end{definition}

\begin{theorem} \label{lim-P}With respect to weak convergence on $\cP(\M)$,
\begin{align*}
\PP^\beta \to \ \PP^0\quad&\text{as }\beta\to0,\\
\PP^\beta \to \ \PP^\infty\quad&\text{as }\beta\to\infty.
\end{align*}
\end{theorem}

\begin{proof} It is easy to see and well-known (cf.~\cite[Ex.~3]{Lynch}) that 
$$\Q^\beta \to \ \PP^0\quad\text{as }\beta\to0, \qquad
\Q^\beta \to \ \PP^\infty\quad\text{as }\beta\to\infty.$$
Indeed, the first assertion follows from  \eqref{particle} and the second one from \eqref{Ferguson}.
Thanks to the continuity of the conjugation map $\Conj:=\Conj_\cP$, this yields
\begin{align*}
\lim_{\beta\to\infty}\PP^\beta=
\lim_{\beta\to\infty}\Conj_*\Q^\beta=\Conj_*\Big(\lim_{\beta\to\infty}\Q^\beta\Big)=\Conj_*\delta_\mm=
\delta_{\Conj(\mm)}=\delta_\mm=\PP^\infty.
\end{align*}
Similarly,
\begin{align*}
\lim_{\beta\to0}\PP^\beta=
\lim_{\beta\to0}\Conj_*\Q^\beta=\Conj_*\Big(\lim_{\beta\to0}\Q^\beta\Big)=\Conj_*\PP^0=\PP^0
\end{align*}
where the last equality follows from
$$\Conj_*\delta_x=\delta_{x^\c}$$
according to Proposition \ref{conj-Diracs} (where $x^\c$ denotes the antipodal point of $x$) and the fact that $x^\c$ is distributed according to $\mm$ if $x$ is so, cf.~\cite[Cor.~3.27, Rem.~3.28]{Lorenzo0}. 
\end{proof}

\subsection{Large Deviations}
The ultimate justifications for our heuristic interpretations \eqref{heur} and \eqref{Qbeta} of the Dirichlet Ferguson measure and of the entropic measure 
are given by the large deviation results to be presented now. They provide a precise %meaning for 
assertion on the asymptotic behavior 
%which in simple terms is summarized as
\begin{align*}
\Q^\beta(A) &\sim  e^{-\beta\, \inf_{\nu\in A} \Ent(\mm\mid\nu)}\\
\PP^\beta(A)  &\sim e^{-\beta\,\inf_ {\mu\in A} \Ent(\mu\mid\mm)}
%\end{align*}
%\begin{align*}
%\begin{array}{c}
%\Q^\beta(A) \sim  e^{-\beta\, \underset{\nu\in A}\inf \Ent(\mm\mid\nu)}\\%\qquad\text{as }\beta\to\infty,\\
%\PP^\beta(A)  \sim e^{-\beta\, \underset{\nu\in A}\inf \Ent(\mu\mid\mm)}
%\end{array}
%\qquad\text{as }\beta\to\infty.
\end{align*}
as $\beta\to\infty$.

\begin{theorem} The Dirichlet-Ferguson measures $\Q^\beta$ for $\beta\to\infty$  satisfy a \emph{Large Deviation Principle} with  rate function
\begin{equation*}
I(\nu):=\Ent(\mm\mid\nu).
\end{equation*}
That is, for every closed set $F\subset\cP(\M)$,
\begin{equation}
\limsup_{\beta\to\infty}\frac1\beta\log\Q^\beta(F)\le - I(F),
\end{equation}
and for every open set $G\subset\cP(\M)$,
\begin{equation}
\liminf_{\beta\to\infty}\frac1\beta\log\Q^\beta(G)\ge - I(G).
\end{equation}
Here and in the sequel, $I(A):=\inf_{\nu\in A}I(\nu)$ for each set $A\subset\cP(\M)$.
\end{theorem}
In dimension $n=1$, more precisely, for $\M=[0,1]$, this result was proven in \cite{Lynch}. (The rate function $I(\nu)$ there, however, was erroneously defined  to be $+\infty$ whenever $\nu\not\ll\mm$.)
 Our proof for general $\M$ follows the argumentation in \cite{Ganesh} where a more sophisticated Large Deviation Principle is derived. As pointed out there, it is unclear whether the final results imply each other. Their proof, however, carries over easily to our setting.  
\begin{proof}
Lemma 2 and Lemma 3 as well as the argumentation in the proof of Theorem 1 of  \cite{Ganesh} carry over to our situation where now $\mu_n$ denotes a measure-valued random value with distribution $\Q^n$ or, in the notation of \cite{Ganesh}, ${\mathfrak D}(n\mm)$
--- with $n$ now in the place of $\beta$ (as usual in our paper.)

The only significant change required is to Lemma 1 there. The defining equation (9) there has to be replaced by
$$\alpha_n(A_j):=n\, \mm(A_j), \qquad\forall n\in\N, \quad j=1,\ldots,k.$$
As a consequence, 
$$Z_n^j\sim {\mathcal G}(n\,\mm(A_j),1)$$
and
\begin{equation*}\lambda_n^j(\theta)=
\begin{cases}
-n\,\mm(A_j)(1-\theta), \quad & \text{if }\theta<1\\
+\infty, &\text{else.}
\end{cases}
\end{equation*}
Thus trivially (even without passing to the limit),
\begin{equation*}\lambda_j(\theta):=\lim_{n\to\infty}\frac1n\lambda_n^j(\theta)
=
\begin{cases}
-\mm(A_j)(1-\theta), \quad & \text{if }\theta<1\\
+\infty, &\text{else.}
\end{cases}
\end{equation*}
The rest of the argumentation then follows verbatim \cite{Ganesh}.
\end{proof}

\begin{theorem} The entropic measures $\PP^\beta$ for $\beta\to\infty$  satisfy a \emph{Large Deviation Principle} with  rate function
\begin{equation*}
J(\mu):=\Ent(\mu\mid\mm).
\end{equation*}
That is, for every closed set $F\subset\cP(\M)$,
\begin{equation}
\limsup_{\beta\to\infty}\frac1\beta\log\PP^\beta(F)\le - J(F),
\end{equation}
and for every open set $G\subset\cP(\M)$,
\begin{equation}
\liminf_{\beta\to\infty}\frac1\beta\log\PP^\beta(G)\ge - J(G).
\end{equation}
Here, $J(A):=\inf_{\nu\in A}J(\nu)$ for each set $A\subset\cP(\M)$.
\end{theorem}

\begin{proof}
This is an immediate consequence of the previous Theorem together with the continuity of the conjugation map and the identity of Theorem \ref{conj-entr} for the relative entropy of conjugated measures.
\end{proof}

\section{Wasserstein Dirichlet Form and Associated Diffusion}
Our goal now is to introduce a `canonical' \emph{Wasserstein Dirichlet form} $\cE_\W$ and associated \emph{Wasserstein diffusion process} on $\cP(\M)$. We will define the former as the relaxation on $L^2(\cP(\M), \PP^\beta)$  of the pre-Dirichlet form 
\begin{align*}
{\mathcal E}^0_\W(f):=\frac12\int_{{\mathcal P}(\M)}
 \big\|\nabla_\W f\big\|^2(\mu)
 \, d\PP^\beta(\mu),
\end{align*}
defined on {cylinder functions on  $\cP(\M)$},
where $\nabla_\W$ denotes the Wasserstein gradient  in the sense of the Otto calculus and where $\PP^\beta$ is the entropic measure with parameter $\beta>0$.

\begin{definition} \begin{enumerate}[\rm (i)]
\item A function $f:  {\mathcal P}(\M)\to\R$ is called \emph{cylinder function} if
there exist $ k\in\N, F\in\C^1(\R^k)$ and $\vec V=(V_1,\ldots,V_k)\in \C^1(\M,\R^k)$ such that
\begin{equation}\label{cyl}
f(\mu)=F\left(\int_\M \vec V\,d\mu\right).\end{equation}
Here and throughout the sequel $\int_\M\vec V\,d\mu:=\big(\int_\M V_1\,d\mu\,\ldots, \int_\M V_k\,d\mu\big)$.
The set of such cylinder functions is denoted by $\text{\sf Cyl}({\mathcal P}(\M))$.

\item For a cylinder function $f$ given as in \eqref{cyl}, its \emph{Otto-Wasserstein gradient} at the point $\mu\in\cP(\M)$ is given by
$$\nabla_\W f(\mu)=
\sum_{i=1}^k \Big(\partial_i F\Big)\bigg(\int_\M \vec V\,d\mu\bigg)
\nabla V_i, 
$$
(cf. \cite{Ot01, AGS05}) 
and the squared norm of the latter at the point $\mu\in\cP(\M)$  is given by 
$$\big\|\nabla_\W f\big\|^2(\mu):=
\sum_{i,j=1}^k \Big(\partial_i F\cdot \partial_j F\Big)\bigg(\int_\M \vec V\,d\mu\bigg)
\int_\M \langle\nabla V_i,\nabla V_j\rangle\,d\mu
.$$
\end{enumerate}
\end{definition}
\begin{definition}
\begin{enumerate}[\rm (i)]
\item Define a quadratic form ${\mathcal E}^0_\W: L^2({\mathcal P}(\M), \PP^\beta) \to\R$ with domain 
\begin{align*}
\Dom({\mathcal E}^0_\W)
:=\text{\sf Cyl}({\mathcal P}(\M))
\end{align*}
by
\begin{align*}
{\mathcal E}^0_\W(f):=\frac12\int_{{\mathcal P}(\M)}
 \big\|\nabla_\W f\big\|^2(\mu)
 \, d\PP^\beta(\mu).
\end{align*}

\item Denote its relaxation in  $L^2({\mathcal P}(\M), \PP^\beta)$ by ${\mathcal E}_\W$,
that is,
$${\mathcal E}_\W(f)=\liminf_{g\to f \text{ in }L^2} {\mathcal E}^0_\W(g).$$
and the domain of the latter by $ {\mathcal F}_\W$.
\end{enumerate}
\end{definition}

\begin{lemma} \begin{enumerate}[\rm (i)]
\item ${\mathcal E}_\W$ is lower semicontinuous on $L^2({\mathcal P}(\M), \PP^\beta)$.

\item ${\mathcal E}_\W$ is quadratic.

\item $\text{\sf Cyl}({\mathcal P}(\M))$ is ${\mathcal E}_\W$-dense in ${\mathcal F}_\W$.
\end{enumerate}
\end{lemma}

\begin{proof}
(i) Standard argument based on diagonal sequences.

(ii) For given $f,g\in {\mathcal F}_\W$ there exist sequences $(f_n), (g_n)$ in $\Dom({\mathcal E}^0_\W)$ with
$$f_n\to f, \ g_n\to g \ \text{in } L^2 \quad \text{and}\quad {\mathcal E}^0_\W(f_n)\to{\mathcal E}_\W(f), \ {\mathcal E}^0_\W(g_n)\to{\mathcal E}_\W(g).$$
Thus by lower semicontinuity of ${\mathcal E}_\W$ and quadraticity of ${\mathcal E}^0_\W$,
\begin{align*}
{\mathcal E}_\W(f+g)&+{\mathcal E}_\W(f-g)\le \liminf_n\big[ {\mathcal E}^0_\W(f_n+g_n)+{\mathcal E}^0_\W(f_n-g_n)\big]\\
&=2 \liminf_n\big[ {\mathcal E}^0_\W(f_n)+{\mathcal E}^0_\W(g_n)\big]=
2 \big[ {\mathcal E}_\W(f)+{\mathcal E}_\W(g)\big].
\end{align*}
Applying the same argument to $f+g$ and $f-g$ in the place of $f$ and $g$ yields the reverse inequality
${\mathcal E}_\W(f)+{\mathcal E}_\W(g)\le
\frac12 \big[ {\mathcal E}_\W(f+g)+{\mathcal E}_\W(f-g)\big]$.

(iii) For given $f\in {\mathcal F}_\W$ choose $(f_n)$ in $\Dom({\mathcal E}^0_\W)$ with
$f_n\to f$ in $L^2$ and ${\mathcal E}^0_\W(f_n)\to{\mathcal E}_\W(f)$. Then also $\frac{f+f_n}2\to f$ in $L^2$ and thus by lower semicontinuity $\liminf_n {\mathcal E}_\W(\frac{f+f_n}2)\ge{\mathcal E}_\W(f)$. By quadraticity therefore
\begin{align*}
0\le \frac12  \limsup_n {\mathcal E}_\W(f-f_n)\le \limsup_n \Big[{\mathcal E}_\W(f)+{\mathcal E}_\W(f_n)-2{\mathcal E}_\W\Big(\frac{f+f_n}2\Big)\Big]\\
\le{\mathcal E}_\W(f)+ \lim_n {\mathcal E}^0_\W(f_n)-2{\mathcal E}_\W(f)
\le0.
\end{align*}
\end{proof}

\begin{theorem} $({\mathcal E}_\W, {\mathcal F}_\W)$ is a regular
Dirichlet form.
\end{theorem}

\begin{proof}
One easily verifies that $\cE_\W(f\wedge 1)\le \cE_\W(f)$ for all $f$. Thus the claim follows from the previous Lemma.
\end{proof}

It is easy to see that $({\mathcal E}_\W, {\mathcal F}_\W)$ equally well can be defined as the relaxation of the pre-Dirichlet form $\cE_\W^0$ with domain restricted to the $\C^\infty$-\emph{cylinder functions}, that is, $f$ with $f(\mu)=F\left(\int_\M \vec V\,d\mu\right)$
for 
$F\in\C^\infty(\R^k)$ and $\vec V=(V_1,\ldots,V_k)\in \C^\infty(\M,\R^k)$.

Let us compare this energy functional $\cE_\W$ with the so-called Cheeger energy.

\begin{definition}
\begin{enumerate}[\rm (i)]
\item Define a convex functional ${\mathcal E}^0_{\sf C\!h}: L^2({\mathcal P}(\M), \PP^\beta) \to\R$ with domain 
$\Dom({\mathcal E}^0_\Ch):=\text{\sf Lip}({\mathcal P}(\M), \W_2)$
by
\begin{align*}
{\mathcal E}^0_\Ch(f):=\frac12\int_{{\mathcal P}(\M)} |{\sf lip}_\W f|^2(\mu)
\, d\PP^\beta(\mu)
\end{align*}
for $f\in \Dom({\mathcal E}^0_\Ch)$ and ${\mathcal E}^0_\Ch(f):=\infty$ else,
where 
$${\sf lip}_\W f(\mu):=\limsup_{\rho,\nu\to \mu, \ \rho\not=\nu}\frac{|f(\rho)-f(\nu)|}{\W_2(\rho,\nu)|}$$ denotes the asymptotic Lipschitz constant of $f$.

\item The Cheeger energy ${\mathcal E}_\Ch$ is the relaxation of ${\mathcal E}^0_\Ch$ in  $L^2({\mathcal P}(\M), \PP^\beta)$, that is,
$${\mathcal E}_\Ch(f)=\liminf_{g\to f \text{ in }L^2} {\mathcal E}^0_\Ch(g).$$
\end{enumerate}
\end{definition}

\begin{lemma}[{\cite[Prop.~4.9]{Savare}}]  $\text{\sf Cyl}({\mathcal P}(\M))\subset{\sf Lip}({\mathcal P}(\M))$, and $$\|\nabla_\W f\|(\mu)={\sf lip}_\W f(\mu)$$
for all $f\in \text{\sf Cyl}({\mathcal P}(\M))$ and all $\mu\in\cP(\M)$.
\end{lemma}

\begin{theorem}[{\cite[Thm.~6.2]{Savare}}]  \label{e=e}
 On $L^2({\mathcal P}(\M), \PP^\beta)$,
\begin{equation}\label{domina}{\mathcal E}_\Ch = {\mathcal E}_\W.\end{equation}
\end{theorem}
In particular, thus the metric measure space 
$(\cP(\M), {\sf W}_2, \PP^\beta)$ is infinitesimally Hilbertian.

\begin{corollary}[{\cite[Prop.~4.11, 4.14]{AGS14a}}] The regular
Dirichlet form
 $({\mathcal E}_\W, {\mathcal F}_\W)$ 
 \begin{itemize}
 \item is strongly local;
\item admits a  carr\'e du champ $\Gamma_\W$, more precisely, 
$$\Gamma_\W(f)=|{\sf lip}_\W f|_*^2= \|\nabla_\W f\|^2_*$$
(here the ${}_*$ refers to the $L^2$-relaxation of the respective quantities)
for $f\in\cF_\W$, and $\Gamma_\W(f)\le |{\sf lip}_\W f|^2$ for  $f\in \text{\sf Lip}({\mathcal P}(\M), \W_2)$;

\item has intrinsic metric 
${\sf d}_\W\ge \W_2$.
\end{itemize}
\end{corollary}

\begin{corollary} There exists  a strong Markov process with a.s.~continuous trajectories (`diffusion process')  properly associated with $({\mathcal E}_\W, {\mathcal F}_\W)$.
\end{corollary}
\begin{definition}
The process $\big((\rho_t)_{t\ge0}, ({\mathbf P}_\mu)_{\mu\in\cP(\M)}\big)$ properly associated with $({\mathcal E}_\W, {\mathcal F}_\W)$
 will be called \emph{Wasserstein diffusion.}
\end{definition}

A crucial question remains open: 
\begin{open}  Is $\mathcal E_\W=\cE^0_\W$ on $\Dom(\cE^0_\W)$? 

Or, in other words, is $(\cE^0_\W,\Dom(\cE^0_\W))$  closable?
\end{open}
 In dimension 1, an affirmative answer will be discussed in the paragraph below.
In higher dimensions one would at least want to know whether the energy is not identically zero.
Indeed, in many `singular' cases the relaxation procedure produces vanishing energy functionals, e.g. 
\begin{itemize}
\item $\cE\equiv0$ for $\cE(f):=\liminf_{g\to f} \cE^0(g)$ if $\cE^0(f):=\int_\R |f'|^2d{\sf p}$ with $\Dom(\cE^0)=\C^1(\R)$ and ${\sf p}\in\cP(\R)$ has no absolutely continuous part;

the same holds true for $\cE^0(f):=\int_\R |\text{lip} f |^2d{\sf p}$ with $\Dom(\cE^0)=\Lip(\R)$;
\item $\cE\equiv0$ for $\cE(f):=\liminf_{g\to f} \cE^0(g)$ if $\cE^0(f):=\int_{S} |\text{lip} f|^2d{\sf p}$ where $S\subset\R^2$ is the Sierpinski gasket, $\Dom(\cE^0)=\Lip(S)$ and $\sf p$ is the Hausdorff measure on $S$.
\end{itemize} 

\begin{open} Is $\mathcal E_\W\not\equiv 0$ on its domain? 

Or, equivalently, does the 
Wasserstein diffusion
$\big((\rho_t)_{t\ge0}, ({\mathbf P}_\mu)_{\mu\in\cP(\M)}\big)$ satisfy
$$\PP^\beta\Big\{\mu: \ {\mathbf P}_\mu\big\{\exists t>0: \rho_t\not= \rho_0\big\}>0\Big\}>0 \ ?$$
\end{open}
An affirmative answer to this will be the main result of the final Section 5.

\paragraph{Wasserstein Diffusion on One-dimensional Spaces}

Let us consider the case $\M=\mathbb S^1$ or $\M=[0,1]$ as treated in \cite{RS09} (in particular, Thm.~7.17, 7.25).

\begin{theorem}[Change-of-variable formula] The entropic measure $\PP^\beta$ is {\em quasi-invariant}  under push-forwards $\mu\mapsto h_*\mu$ by means of smooth diffeomorphisms $h$: 
$$\frac{d\PP^\beta(h_*\mu)}{d\PP^\beta(\mu)}=X_h^\beta\cdot Y_h(\mu)$$
where the density consists of of two terms. The first one
$$X_h^\beta (\mu)=\exp\left(\beta\int_{\M}\log h'(t)d\mu(t)\right)$$
can be interpreted as $\exp(-\beta\Ent(h_*\mu|\mm)) / \exp(-\beta\Ent(\mu|\mm))$ in accordance with our formal interpretation
(\ref{heur}).
The second one
$$Y_h(\mu)=
\prod_{I\in{\sf gaps}(\mu)}\frac{\sqrt{h'(I_-)\cdot h'(I_+)}}{|h(I)|/|I|}
$$
takes into account the action of $h$ on 
%can be interpreted as the change of variable formula for the (non-existing) measure $\PP^*$.
%Here 
the set ${\sf gaps}(\mu)$  consisting of the intervals $I=(I_-,I_+)\subset \M$ of maximal length with $\mu(I)=0$.
\end{theorem}
For an alternative proof for this change-of-variable formula, see also \cite{RYZ08}.
Such a change-of-variable formula directly leads to an integration-by-parts formula and thus proves closability of the associated pre-Dirichlet form.

\begin{corollary} The pre-Dirichlet form $(\cE^0_\W,\Dom(\cE^0_\W))$ is closable.
\end{corollary}

\begin{theorem} The Dirichlet form $(\cE_\W,\Dom(\cE_\W))$ on one-dimensional spaces has the following properties:
\begin{itemize}
\item Its intrinsic metric  coincides with $\W_2$ \cite[Cor.~7.29]{RS09}.
\item It has  spectral gap $\ge \beta/2$ and satisfies a logarithmic Sobolev inequality, \cite{Stannat}
\item It satisfies no synthetic lower Ricci bound, \cite{Chodosh}.
\item The Wasserstein diffusion admits approximations in terms of interacting particle systems, \cite{AR07}, \cite{Sturm-particle}.
\end{itemize}
\end{theorem}

\paragraph{Related Approaches} 
\begin{itemize}

 \item[(i)] 
 An alternative construction of a reversible diffusion on the Wasserstein space $\cP(\M)$ has been proposed by Dello Schiavo \cite{DelloSchiavo}. This 'Dirichlet Ferguson diffusion' is properly associated with the strongly local, regular Dirichlet form $\cE_\cP$ on $L^2(\cP(\M),\Q^\beta)$
obtained as closure of the
 pre-Dirichlet form 
 \begin{align*}
{\mathcal E}_\cP^0(f):=\frac12\int_{{\mathcal P}(\M)}
 \big\|\nabla_\W f\big\|^2(\mu)
 \, d\Q^\beta(\mu),
\end{align*}
where $\nabla_\W$ again denotes the Wasserstein gradient  in the sense of the Otto calculus and where $\Q^\beta$ is the Dirichlet-Ferguson measure with parameter $\beta>0$.
  
\item[(ii)]  A closely related construction of a reversible diffusion on the configuration space $\Upsilon(\M)\subset {\mathcal M}(\M)$ has been presented by Albeverio, Kondratiev, R\"ockner \cite{AKR98}.
This diffusion is properly associated with the strongly local, regular Dirichlet form $\cE_\Upsilon$ on $L^2(\cP(\M),\pi_{\lambda,\mm})$
obtained as closure of the
 pre-Dirichlet form 
 \begin{align*}
{\mathcal E}_\Upsilon^0(f):=\frac12\int_{{\Upsilon}(\M)}
 \big\|\nabla_\W f\big\|^2(\mu)
 \, d\pi_{\lambda,\mm}(\mu),
\end{align*}
where $\nabla_\W$ again denotes the Wasserstein gradient  in the sense of the Otto calculus and where $\pi_{\lambda,\mm}=\int_0^\infty \pi_{s\mm}d\lambda(s)$ is the mixture w.r.t.~$\lambda\in\cP(\R_+)$ of Poisson measures $\pi_{s\mm}$ with intensity $s\mm$.

 \item[(iii)]  As modifications of the Wasserstein diffusion of \cite{RS09} on 1-dimensional spaces, Konarovskyi and von Renesse \cite{KR17,KR18} present and analyze the 
modified massive Arratia flow  and  the
coalescing-fragmentating Wasserstein dynamics.

\item[(iv)]  Chow and Gangbo \cite{Gangbo} study a stochastic process on $\cP(\M)$ `modeled after Brownian motion' and generated by a `partial Laplacian'.

 \end{itemize}

\section{Being in Motion}

The main challenge which remains now in the multidimensional case is to prove that $\cE_\W\not=0$ or, in other words, that the associated stochastic processes really moves.

We will prove that $\cE_\W(f)\not=0$ for all non-vanishing $f$ in two major classes of functions on $\M$. 
Our argument is based on the existence of suitable families of isometries which is granted 
for the $n$-sphere and for the $n$-torus.

\subsection{Lipschitz Families of Isometries}

\begin{definition} (i) Given a Riemannian compactum $\M$, a family $(\Phi_t)_{t\in[0,1]}$ is called \emph{Lipschitz family of isometries} iff
\begin{itemize}
\item $\Phi_0={\sf Id}$  
\item $\Phi_t:\M\to\M$ is an isometry for every $t\in [0,1]$
\item $\d(\Phi_s(y), \Phi_t(y))\le L\cdot  |t-s|$ for all $s,t\in [0,1]$, all $y\in\M$, and some $L\in\R$.
\end{itemize}

(ii) A Riemannian compactum $\M$ is called \emph{Lipschitz transitive} iff for any pair of points $x_0,x_1\in\M$,  %with $d(x_0,x_1)< {\sf inf}(\M)$, 
there exists a Lipschitz family of isometries 
$(\Phi_t)_{t\in[0,1]}$ with $L=\d(x_0,x_1)$ and $\Phi_1(x_0)=x_1$.
\end{definition}

\begin{lemma}\label{iso} 
\begin{enumerate}[\rm (i)]
\item The $n$-torus is Lipschitz transitive.

\item The $n$-sphere is Lipschitz transitive.
\end{enumerate}
\end{lemma}

\begin{proof} (i) In the case of the torus, the additive structure of $\R^n$ modulo $\Z^n$ allows  easily to verify that the following maps will do the job:
$$\Phi_r(y)=y+r(x_1-x_0).$$

(ii) Now consider the case of the sphere, regarded as $$\M=\bigg\{y=\big(y^{(1)}, \ldots, y^{(n+1)}\big) \in \R^{n+1}: \sum_{i=1}^{n+1} \big|y^{(i)} \big|^2=1\bigg\}.$$ Without restriction, we may assume that $x_0$ and $x_1$ lie in the $e_1$-$e_2$-plane, say,
$$x_0=(1,0,\ldots, 0), \qquad x_1=(\cos\varphi, \sin\varphi, 0,\ldots,0).$$
Then we define the $\Phi_t$'s as rotations of the $e_1$-$e_2$-plane. More precisely,
$$\Phi_t: \big(r \cos\psi, r\sin\psi, y_3, \ldots,y_{n+1}\big) \mapsto \big(r \cos(\psi+t\varphi), r\sin(\psi+ t \varphi), y_3, \ldots,y_{n+1}\big) .$$
Obviously, this defines isometries of the sphere with the requested properties.
\end{proof}

\begin{definition}
A curve $(\mu_t)_{[0,1]}$ in the metric space  $(\cP(\M), {\sf W}_2)$ is called 2-absolutely continuous, briefly $(\mu_t)_{[0,1]}\in{\mathcal AC}_2\big([0,1]; \big(\cP(\M), \sf W_2\big)\big)$, 
iff there exists $h\in L^2([0,1],\R)$ with
$$
{\sf W}_2(\mu_s,\mu_t)\le \int_s^t h(r)dr \qquad \forall s,t\in[0,1], s<t.$$

\end{definition}
\begin{lemma}\label{ac2-compr}
Let a Lipschitz family $(\Phi_t)_{t\in [0,1]}$ of isometries be given. Then
\begin{itemize}
\item[\rm (i)] for each $t$, the measure $\PP^\beta$ is invariant under the map $\hat \Phi_t: \mu\mapsto (\Phi_t)_*\mu$, i.e.
$$(\hat\Phi_t)_*\PP^\beta=\PP^\beta.$$
\item[\rm (ii)] The map $\hat\Phi: \mu \mapsto \big( (\Phi_t)_*\mu\big)_{[0,1]}$ maps $\cP(\M)$ into $\cP\Big({\mathcal AC}_2\Big([0,1]; \big(\cP(\M), \sf W_2\big)\Big)\Big)$. Thus
\begin{equation}\label{P-vec}
\vec\PP^\beta:=\hat\Phi_*\PP^\beta\in \cP\Big({\mathcal AC}_2\Big([0,1]; \big(\cP(\M), \sf W_2\big)\Big)\Big).
\end{equation}
Moreover, for every $t$,
$$(e_t)_*\vec\PP^\beta=(\hat\Phi_t)_*\PP^\beta=\PP^\beta.$$
\end{itemize}

\end{lemma}

\begin{proof}
(i) follows from $\hat\Phi_t$-invariance of the measure $\Q^\beta$ and $\hat\Phi_t$-invariance of the conjugation map $\mu\mapsto \mu^\c$ (which in turn follows from 
$\Phi_t$-invariance of $\d$ and $\mm$).

(ii) For all $s,t$
\begin{align*}
{\sf W}_2^2(\mu_s,\mu_t)\le \int_M \d^2(\Phi_s(x),\Phi_t(x))d\mu(x)\le L^2\cdot |s-t|^2.
\end{align*}
\end{proof}

\begin{definition} \label{def-anti}
A function $V:\M\to\R$ is called \emph{antisymmetric} iff 
$$V\circ \Phi_1=-V$$ for a suitable  Lipschitz family $(\Phi_t)_{t\in [0,1]}$ of isometries,
and 
 a function $f:\cP(\M)\to\R$ is called \emph{antisymmetric} iff 
$$f\circ \hat\Phi_1=-f.$$
\end{definition}

\begin{remark}
\begin{enumerate}[\rm (i)]
\item If $V:\M\to\R$ is antisymmetric then $f:\cP(\M)\to\R$ defined by $f(\mu):=\int_\M V\,d\mu$ is antisymmetric.

\item More generally, if $V_i:\M\to\R$ are antisymmetric for $i=1,\ldots,k$ (w.r.t.~the same Lipschitz family $(\Phi_t)_{t\in [0,1]}$ of isometries)
and $F:\R^k\to \R^k$ is antisymmetric in the usual sense, i.e.~$F(-x)=-x$, then $f:\cP(\M)\to\R$ defined by $f(\mu)=F\big(\int_\M V_1\,d\mu\,\ldots, \int_\M V_k\,d\mu\big)$
is antisymmetric.
\end{enumerate}
\end{remark}

\subsection{Two Lower Bounds on the Wasserstein Energy}

\begin{definition} 
Let $f$ be a Borel function on  $\cP(\M)$. Then a Borel function $g$ on $\cP(\M)$ is called \emph{weak upper gradient} for $f$ iff for every 
$$\mathbf P\in \cP\Big({\mathcal AC}_2\Big([0,1]; \big(\cP(\M), \sf W_2\big)\Big)\Big) \quad\text{with}\quad\sup_t \ (e_t)_*\mathbf P\le C\cdot  \PP^\beta$$
and for $\mathbf P$-a.e.~$(\mu_t)_{[0,1]}\in{\mathcal AC}_2\Big([0,1]; \big(\cP(\M), \sf W_2\big)\Big)$
\begin{align*}
\Big| f(\mu_1)- f(\mu_0)\Big|\le \int_0^1 g(\mu_t)\cdot \big|\dot \mu_t\big|\,dt.
 \end{align*}
\end{definition}

\begin{proposition}\label{key-est} Let $(\Phi_t)_{t\in [0,1]}$ be a Lipschitz family of isometries. Then for every Borel function $f$ on   $\cP(\M)$,
\begin{equation}
\cE_\W(f)\ge \frac1{2L^2}\int_{\cP(\M)}\Big| 
f\big((\Phi_1)_*\mu\big)- f\big(\mu\big)
\Big|^2d\PP^\beta(\mu).
\end{equation}
\end{proposition}

\begin{proof} By definition and the fact that $|\dot \mu_t|\le  L$, for every weak upper gradient $g$ and for $\mathbf P=\vec\PP^\beta$ as defined in \eqref{P-vec},
\begin{align*}
\Big| f(\mu_1)- f(\mu_0)\Big|\le L\int_0^1 g(\mu_t)\,dt
 \end{align*}
 for
$\vec\PP^\beta$-a.e.~$(\mu_t)_{t\in[0,1]}\in{\mathcal AC}_2\big([0,1]; \big(\cP(\M), \sf W_2\big)\big)$, \cite[Def.~5.4 and Thm.~6.2]{AGS14}. Integrating the squared inequality w.r.t.~$\vec\PP^\beta$ and using the fact that 
$(e_t)_*\vec\PP^\beta=\PP^\beta$
yields
\begin{align*}
\int_{{\mathcal AC}_2}\Big| f(\mu_1)- f(\mu_0)\Big|^2 d\vec\PP^\beta\big((\mu_t)_t\big)&\le
L^2 \int_{{\mathcal AC}_2}  \bigg|\int_0^1 g(\mu_t)\,dt\bigg|^2d\vec\PP^\beta\big((\mu_t)_t\big)\\
&\le L^2\int_{\cP(\M)} g^2(\mu)d\PP^\beta(\mu).
\end{align*}
Taking the infimum w.r.t.~$g$ thus yields the claim
\begin{align*}
\cE_\Ch(f)&\ge \frac1{2L^2}\int_{{\mathcal AC}_2}\Big| f(\mu_1)- f(\mu_0)\Big|^2 d\vec\PP^\beta\big((\mu_t)_t\big)\\
&= \frac1{2L^2}\int_{\cP}\Big| f((\Phi_1)_*\mu)- f(\mu)\Big|^2 d\PP^\beta(\mu).
\end{align*}
\end{proof}

\begin{theorem}\label{nix} For every antisymmetric Borel function $f$ on $\cP(\M)$,
\begin{equation}\label{anti-symm}
\cE_\W(f)\ge \frac2{L^2}\int_{\cP(\M)} f^2\big(\mu\big)
\, d\PP^\beta(\mu).
\end{equation}
In particular,
\begin{align*}\cE_\W(f)=0\quad
&\Longleftrightarrow\quad f=0\ \PP\text{-a.s.}\quad
\Longleftrightarrow\quad f\equiv0.
\end{align*}
\end{theorem}

\begin{proof} By symmetry, $|f((\Phi_1)_*\mu)- f(\mu)|^2=4 f^2(\mu)$. With this,  the 
first assertion follows readily from the previous Proposition \ref{key-est}. 

For the second assertion, observe that 
 \begin{itemize}
\item $\cE_\W(f)=0\  \Rightarrow \  f=0\ \PP$-a.s. \quad follows from \eqref{anti-symm} 
\item $ f=0\ \PP^\beta\text{-a.s.}\ \Rightarrow\  f\equiv0$ \quad follows from Proposition \ref{full-supp}
\item $f\equiv0\ \Rightarrow \ \cE_\W(f)=0$ \quad is trivial.
\end{itemize}
\end{proof}

\begin{theorem} Assume that $\M$ is Lipschitz transitive. Then for every non-constant Lipschitz function $V:\M\to\R$, the function $f:{\mathcal P}(\M)\to\R$,
$  \mu\mapsto  \int_\M V\,d\mu$ has nonvanishing Wasserstein energy $\mathcal E_\W(f)$.
\end{theorem}
More explicitly, for every $x_0,x_1\in\M$ and $\epsilon>0$ with $|V(x_1)-V(x_0)| \ge {\sf Lip}V\cdot \Big[\frac12d(x_0,x_1)+2\epsilon\Big]$ and with
$\eta:=\PP^\beta\big( {\mathbb B}_\epsilon(\delta_x)\big)>0$,
\begin{equation}
\mathcal E_\W(f)\ge \frac\eta8 ({\sf Lip}V)^2.
\end{equation}

\begin{proof}  Obviously, $\eta$ as defined above is independent of $x$. The positivity of it was proven in Theorem \ref{full-supp}.

Let points $x_0,x_1$ be given as above and let $(\Phi_t)_t$ be the associated family of isometries with $\Phi_1(x_0)=x_1$. Then
$$\mu\in {\mathbb B}_\epsilon(\delta_{x_0}) \ \Leftrightarrow \ (\Phi_1)_*\mu\in {\mathbb B}_\epsilon(\delta_{x_1}).$$
Furthermore, if $\mu\in {\mathbb B}_\epsilon(\delta_{x_0})$ and $\nu\in {\mathbb B}_\epsilon(\delta_{x_1})$, then by triangle inequality and Kantorovich duality
\begin{align*}\left|\int V d\nu-\int Vd\mu\right|&\ge 
\Big|V(x_1)-V(x_0)\Big|-\Big|\int V d\nu- V(x_1)\Big|-\Big|\int V d\mu- V(x_0)\Big|\\
&\ge \Big|V(x_1)-V(x_0)\Big|- {\sf Lip}V\cdot W_1(\nu,\delta_{x_1})-{\sf Lip}V\cdot W_1(\mu,\delta_{x_0})\\
&\ge \Big|V(x_1)-V(x_0)\Big|- 2\epsilon\cdot{\sf Lip}V\\
&\ge\frac12 {\sf Lip}V\cdot d(x_0,x_1)
\end{align*}
where the last inequality is due to  our choice of $x_0,x_1$ and $\epsilon$.
Thus combined with Proposition \ref{key-est}, with $L=d(x_0,x_1)$,
\begin{align*}
\cE_\W(u)&\ge \frac1{2L^2}\int_{\cP(\M)}\Big| 
u\big((\Phi_1)_*\mu\big)- u\big(\mu\big)
\Big|^2d\PP(\mu)\\
&\ge \frac1{2L^2}\int_{{\mathbb B}_\epsilon(\delta_{x_0})}\Big| 
u\big((\Phi_1)_*\mu\big)- u\big(\mu\big)
\Big|^2d\PP(\mu)\\
&\ge \frac1{2 L^2}\Big(\frac12{\sf Lip}V\cdot d(x_0,x_1)\Big)^2\cdot \PP\big( {\mathbb B}_\epsilon(\delta_{x_0})\}\big)\\
&=\frac\eta{8}({\sf Lip}V)^2.
\end{align*}
\end{proof}

\subsection{An Asymptotically Sharp Lower Estimate for the Wasserstein Energy}
Let us finally consider a particularly easy example where the Wasserstein energy indeed is arbitrarily close to the corresponding pre-energy.

\begin{theorem}
Let $\M=\R^n/\Z^n$ be the $n$-torus,  identified with the set $(-1/2,1/2]^n\subset\R^n$, and consider
$f:\cP(\M)\to\R$ defined as $f(\mu)=\int V\,d\mu$ with 
$V(x)=|x_1|$.
Then for every $\epsilon>0$ there exists $\beta_*>0$ %and $\beta^*<\infty$ 
such that for each $\beta\in (0,\beta_*)$,
$$\cE_\W(f)\ge (1-\epsilon)\, \cE_\W^0(f)$$
%where $\cE_\W^0(f)=\frac12\int |\nabla V|^2\,d\mm =\frac12$ and 
where $\cE_\W$ denotes the relaxation of $\cE_\W^0$ in $L^2(\cP(\M),\PP^\beta)$.
\end{theorem}

\begin{proof} Put $\M_\epsilon:=\{x\in\M: 0<x_1<1/2-\epsilon/2\}$ and $\mathcal M_\epsilon:=\{\mu\in\cP(\M): \mu(\M_\epsilon)>1-\epsilon\}$. Then $\mathcal M_\epsilon$ is open in $\cP(\M)$.  Therefore, Theorem \ref{lim-P} implies
\begin{align*}
\liminf_{\beta\to0}\PP^\beta(\mathcal M_\epsilon)\ge \PP^0(\mathcal M_\epsilon)=\mm\big\{x\in\M: \, \delta_x(\M_\epsilon)>1-\epsilon\big\}=\mm(\M_\epsilon)=\frac12-\frac\epsilon2.
\end{align*}
Thus there exists  $\beta_*>0$ such that 
$$\PP^\beta(\mathcal M_\epsilon)\ge \frac12-\epsilon$$
for all $\beta\in (0,\beta_*)$.
For fixed $\beta\in (0,\beta_*)$ we define 
$\PP^\beta_\epsilon(.):=\frac1{C_\epsilon}\, \PP^\beta(\M_\epsilon \cap .)$
with $C_\epsilon:=\PP^\beta(\M_\epsilon)$, and similarly as in \eqref{P-vec},
\begin{equation*}
\vec\PP^\beta_\epsilon:=\hat\Phi_*\PP^\beta_\epsilon\in \cP\Big({\mathcal AC}_2\Big([0,1]; \big(\cP(\M), \sf W_2\big)\Big)\Big)
\end{equation*}
where $\Phi_t(x):=x+te_1$. From this we deduce
that
\begin{align*}
1-2\epsilon\le\frac1\delta\int_0^\delta \partial_1 V(\Phi_t(x))\,d\mu_0(x)=\frac1\delta
\Big( f(\mu_\delta)- f(\mu_0)\Big)\le \frac1\delta\int_0^\delta g(\mu_t)\,dt
 \end{align*}
 for
$\vec\PP^\beta_\epsilon$-a.e.~$(\mu_t)_{t\in[0,1]}\in{\mathcal AC}_2\big([0,1]; \big(\cP(\M), \sf W_2\big)\big)$ and every $\delta\in(0,\epsilon/2)$. %{\color{red} Use better argumentation! } 
Hence,
%$g(\mu)\ge 1-2\epsilon$ for $\PP^\beta$-a.e.~$\mu\in\mathcal M_\epsilon$.
$$1-2\epsilon\le  \int_{\cP(\M)}\frac1\delta\int_0^\delta g((\Phi_t)_*\mu)\,dt\, d\PP^\beta_\epsilon(\mu)=\frac1{C_\epsilon}
 \int_{\mathcal M_\epsilon} g(\mu)\, d\PP^\beta(\mu)
$$
Essentially the same argumentation applies to $\M'_\epsilon:=\{x\in\M: -1/2<x<-\epsilon/2\}$ and $\mathcal M'_\epsilon:=\{\mu\in\cP(\M): \mu(\M'_\epsilon)>1-\epsilon\}$. Obviously,  $\mathcal M_\epsilon$ and $\mathcal M'_\epsilon$ are disjoint, and $\PP^\beta(\mathcal M'_\epsilon)\ge \frac12-\epsilon$. Moreover, $g(\mu)\ge 1-2\epsilon$ for $\PP^\beta$-a.e.~$\mu\in\mathcal M'_\epsilon$.
Therefore,
\begin{align*}
\int_{\cP(\M)}g\,d\PP^\beta&\ge \int_{\mathcal M_\epsilon}g\,d\PP^\beta+ \int_{\mathcal M'_\epsilon}g\,d\PP^\beta\ge (1-2\epsilon)\, \big(
\PP^\beta(\mathcal M_\epsilon)+\PP^\beta(\mathcal M'_\epsilon)\big)\ge(1-2\epsilon)^2.
\end{align*}
Taking the minimal $g$ thus yields
$$\cE_\W(f)=\frac12\int_{\cP(\M)}g^2\,d\PP^\beta\ge\frac12\left( \int_{\cP(\M)}g\,d\PP^\beta\right)^2
\ge \frac12 (1-2\epsilon)^4.$$
On the other hand, obviously
$\cE^0_\W(f)=\frac12\int |\nabla V|^2\,d\mm =\frac12$.
\end{proof}

\paragraph{Acknowledgement.} The author gratefully acknowledges  funding by the Deutsche Forschungsgemeinschaft through the Hausdorff 
Center for Mathematics and through CRC 1060 as well as through SPP 2265. He also thanks Lorenzo Dello Schiavo for providing the simulations in Figure 1.
%\end{acknowledgement}

\end{document}